\documentclass[10pt,reqno]{amsart}

\usepackage{amsmath,amssymb,amsthm,mathtools}
\usepackage[marginratio=1:1,height=600pt,width=470pt,tmargin=100pt]{geometry}
\usepackage{enumitem}
\usepackage{hyperref}
\usepackage{xcolor}
\usepackage{tikz-cd}
\usetikzlibrary{arrows}

\newtheorem{theorem}{Theorem}[section]
\newtheorem{lemma}[theorem]{Lemma}
\newtheorem{proposition}[theorem]{Proposition}
\newtheorem{corollary}[theorem]{Corollary}
\newtheorem{conjecture}[theorem]{Conjecture}
\theoremstyle{definition}

\newtheorem{condition}[theorem]{Condition}
\theoremstyle{remark}
\newtheorem{remark}[theorem]{Remark}
\newtheorem{example}[theorem]{Example}
\numberwithin{equation}{section}

\def\H{\mathrm{H}}
\def\P{\mathrm{P}}

\def\Q{\mathbb{Q}}
\def\C{\mathbb{C}}
\def\A{\mathbb{A}}
\def\N{\mathbb{N}}
\def\GM{{{\mathbb{G}_m}}}
\def\BGM{\mathrm{B}\mathbb{G}_m}
\def\DQ{\mathbb{D}\mathbb{Q}}
\def\BM{\mathrm{BM}}
\def\rd{\mathrm{rd}}
\def\cl{\mathrm{cl}}
\def\ss{\mathrm{ss}}
\def\gr{\mathrm{gr}}
\def\taut{\mathrm{taut}}
\def\ch{\mathrm{ch}}
\def\cBPS{\mathcal{BPS}}
\def\BPS{\mathrm{BPS}}
\def\Pure{\mathrm{Pure}}
\def\Spec{\mathrm{Spec}}
\def\calH{\mathcal{H}}
\def\calU{\mathcal{U}}
\def\calA{\mathcal{A}}
\def\calO{\mathcal{O}}
\def\calE{\mathcal{E}}
\def\calM{\mathcal{M}}
\def\calT{\mathcal{T}}
\def\bbT{\mathbb{T}}
\def\bbS{\mathbb{S}}

\def\g{\mathfrak{g}}
\def\bfH{\mathbf{H}}
\DeclareMathOperator{\vdim}{vdim}

\DeclareMathOperator{\IC}{IC}
\def\Irr{\mathrm{Irr}}
\DeclareMathOperator{\Sym}{Sym}
\let\Im\relax\DeclareMathOperator{\Im}{Im}
\def\Q{\mathbb{Q}}
\def\C{\mathbb{C}}
\def\N{\mathbb{N}}
\def\Z{\mathbb{Z}}
\def\A{\mathbb{A}}
\def\L{\mathbb{L}}
\def\T{\mathbb{T}}

\def\vd{\mathrm{vdim}}
\def\rk{\mathrm{rk}}
\def\vir{\mathrm{vir}}
\def\Tor{\mathrm{Tor}}
\def\Ob{\mathrm{Ob}}
\def\cofib{\mathrm{cofib}}
\def\H{\mathrm{H}}
\def\BM{\mathrm{BM}}

\def\IC{\mathrm{IC}}
\def\Id{\mathrm{Id}}
\def\rd{\mathrm{rd}}
\def\cl{\mathrm{cl}}
\def\calM{\mathcal{M}}
\def\tcalM{\widetilde{\mathcal{M}}}
\def\g{\mathfrak{g}}
\def\n{\mathfrak{n}}
\def\bfA{\mathbf{A}}
\def\U{\mathrm{U}}
\def\BGM{\mathrm{B}\mathbb{G}_m}
\def\B{\mathrm{B}}
\def\GM{\mathbb{G}_m}
\def\Rep{\mathrm{Rep}}
\def\tr{\mathrm{tr}}
\def\DQ{\mathbb{D}\Q}
\def\cBPS{\mathcal{BPS}}
\def\BPS{\mathrm{BPS}}

\def\scal{\mathrm{scal}}
\def\snilp{\mathrm{snil}}
\def\DT{\mathrm{DT}}
\def\Filt{\mathfrak{Filt}}
\def\ev{\mathrm{ev}}

\def\H{\mathrm{H}}
\def\Q{\mathbb{Q}}
\def\C{\mathbb{C}}
\def\A{\mathbb{A}}
\def\N{\mathbb{N}}
\def\L{\mathbb{L}}
\def\T{\mathbb{T}}
\def\GM{\mathbb{G}_m}
\def\DQ{\mathbb{D}\mathbb{Q}}
\def\BM{\mathrm{BM}}
\def\red{\mathrm{red}}
\def\rd{\mathrm{rd}}
\def\nil{\mathrm{nil}}
\def\snil{\mathrm{snil}}
\def\scal{\mathrm{scal}}
\def\taut{\mathrm{taut}}
\def\BPS{\mathrm{BPS}}
\def\BM{\mathrm{BM}}
\def\g{\mathfrak{g}}
\def\gl{\mathfrak{gl}}
\def\U{\mathrm{U}}
\def\ss{\mathrm{ss}}

\def\bfA{\mathbf{A}}
\def\calV{\mathcal{V}}
\def\calE{\mathcal{E}}
\def\ud{{\underline{d}}}
\def\ue{{\underline{e}}}
\def\uf{{\underline{f}}}
\def\Rep{\mathfrak{Rep}}
\def\tRep{\widetilde{\mathfrak{Rep}}}
\def\IQ{IQ}
\DeclareMathOperator{\coker}{coker}
\let\Im\relax\DeclareMathOperator{\Im}{Im}

\title{On the pure part of the cohomology of character varieties}
\author{O. Schiffmann and E. Vasserot}
\date{}

\address{O. Schiffmann: Laboratoire de Math\'ematiques d'Orsay, B\^at 307, Facult\'e des Sciences d'Orsay, Universit\'e Paris-Saclay, 91405 Orsay, France AND Simion Stoilow Institute of Mathematics, P.O. Box 1-764, RO-014700 Bucharest, Romania}
\email{olivier.schiffmann@universite-paris-saclay.fr}

\address{E. Vasserot: IMJ-PRG, UMR7586, Universit\'e Paris-Cit\'e, 75013 Paris, France}
\email{vasserot@imj-prg.fr}

\begin{document}
\maketitle

\begin{abstract}
We prove the purity conjecture of Hausel and Rodriguez-Villegas relating the graded dimension of the pure part of the cohomology of the twisted character variety $M^B(r,d)$ of $GL(r)$ of genus $g$ with the cohomology of a suitable quiver variety associated to the quiver $S_g$ with one vertex and $g$ loops. In particular, this provides an explicit formula for the pure part of the mixed Hodge polynomial of $M^B(r,d)$ in terms of the Kac polynomial $A_{S_g,r}(q)$. Our approach is based on the P=W theorem to transfer the purity conjecture into a statement about the perverse filtration on the cohomology of the stack of Higgs bundles on a genus $g$ curve. We use the theory of BPS Lie algebras, their $\chi$-independence property (proven by Kinjo-Koseki), the recent proof of the stacky version of P=C due to Picovschi and Davison's analytic neighborhood theorem to reduce the conjecture to a comparison problem between (the 'pure' part of) the BPS Lie algebra associated to the stack of Higgs bundles on a genus $g$ curve and the BPS Lie algebra of the stack of representations of the preprojective algebra of the quiver $S_g$. This last step is taken care of by a quiver analog of Markman's tautological generation theorem for the BPS Lie algebra of arbitrary quivers (the latter result is also proven by Davison and Hennecart, see \cite{DH}).
\end{abstract}

\setcounter{tocdepth}{1}
\tableofcontents

%=====================================================================
\medskip\section{Introduction}
%=====================================================================

\medskip\subsection{The purity conjecture}
Let $\Sigma_g$ be a smooth connected and compact topological surface of genus $g\geqslant2$. For $r\geqslant1$ and $d \in (\Z/r\Z)^\times$, the twisted character variety is defined as
\begin{equation}\label{eq:defcharacter}
M^B(r,d)=\Big\{(A_1,B_1,\ldots,A_g,B_g)\in GL_r^{2g}(\C)\;;\;\prod_i[A_i,B_i]=e^{2i\pi\frac{d}{r}}\Big\}/\!\!/GL_r .
\end{equation}
It may be interpreted as the moduli space of $GL_r(\C)$-representations of the once-punctured $\Sigma_g$, with monodromy around the puncture given by $e^{2i\pi \frac{d}{r}}$. The study of the cohomology of this affine algebraic variety was initiated by T. Hausel and F. 
Rodriguez-Villegas in the seminal paper \cite{HRV08}, who showed in particular that $M^B(r,d)$ is smooth, connected, of dimension $2N=2(g-1)r^2+2$ and carries a natural symplectic form. In particular, the mixed Hodge structure on $\H^m(M^B(r,d),\Q)$ has weights in the interval $[m,2m]$. It is of Hodge-Tate type by \cite{S17}. 
The {pure part} of $\H^*(M^B(r,d))$ is 
$$\H^{*}_{\text{pure}}(M^B(r,d),\Q)=\bigoplus_m W_m\H^m(M^B(r,d),\Q)$$ 
It is concentrated in even degrees, since $W_{2k+1}=W_{2k}$ and $W_{2k}\H^{2k+1}=0$.

Let $S_g$ stand for the quiver with one vertex and $g$ loops, and let $A_{S_g,r}(q)\in\N[q]$ 
be the Kac polynomial of $S_g$ counting absolutely indecomposable $S_g$-representations over $\mathbb{F}_q$ of dimension $r$. 
It is of degree $N=(g-1)r^2+1$.

As a special case of their conjectural expression for the full mixed Hodge polynomial of $M^B(r,d)$, Hausel and Rodriguez-Villegas formulated the following {purity conjecture}.

\begin{conjecture}[{\cite{HRV08}}]\label{conj:HRV}
	For all $r\geqslant1$ and $d$ coprime to $r$ we have
	$$\sum_{k}\dim W_{2k}\H^{2k}(M^B(r,d),\Q)\,q^{k}=q^{(g-1)r^2+1}A_{S_g,r}(q^{-1}).$$
	Equivalently, by Poincar\'e duality, the graded dimension of the pure part of compactly supported cohomology is 
	$q^{(g-1)r^2+1}A_{S_g,r}(q)$.
\end{conjecture}

\begin{example}\hfill
	\begin{enumerate}[label=(\alph*),leftmargin=8mm]
		\item
		For $r=1$: $M^B(1,0)=(\C^*)^{2g}$ and the pure part is $\H^0$.
		Indeed $q^{g}A_1^{S_g}(q^{-1})=q^gq^{-g}=1$. 
		\item
		For $r=2$: $A_2^{S_g}(q)=q^{2g-1}(1+q^2+\cdots+q^{2g-2})$, 
		so the conjecture predicts one pure class in each of the degrees $0,4,\ldots,4g-4$ (see \cite{HauselThaddeus2} and \cite{Hua} respectively).
	\end{enumerate}
\end{example}

\medskip

The main result of this paper is a proof of the above purity conjecture.

\begin{theorem}\label{thm:purityconj}
	Conjecture~\ref{conj:HRV} is true.
\end{theorem}

\medskip\subsection{Main ideas of the proof} Let us say a few words about the proof, which combines several recent results about the so-called BPS Lie algebras associated to the stacks of representations of quivers or to the stack of Higgs bundles on smooth complex projective curves.
First of all, we use the P=W theorem proved in \cite{MS} and \cite{HMMS} to translate the problem into the study of the perverse filtration, induced by the Hitchin morphism, of the cohomology of the moduli space $M^D(r,d)$ of semistable Higgs bundles of rank $r$ and degree $d$ 
on a smooth, connected, complex projective curve $C$ of genus $g$. 
More precisely, by Poincar\'e duality the purity conjecture is equivalent to the equality
\begin{equation}\label{eq:intro1}
	\sum_\ell \dim \gr^\P_{N-\ell}\, \H_{2N+2\ell}(M^D(r,d),\Q)q^\ell=A_{S_g,r}(q) ,
\end{equation}
where $P_k$ is the shifted perverse filtration and $\H_*(-,\Q)$ stands for Borel-Moore homology. Next, we recast this in 
terms of the affinized BPS cohomology $\widehat{\mathfrak{g}}^{C}(r,d)$ of $M^{D}(r,d)$.
Note that, for a coprime pair $(r,d)$, the BPS cohomology coincides up to a shift with 
the (co)homology of the stack $\mathfrak{M}^{ss}(r,d)$ of semistable Higgs bundles.
Then we use the $\chi$-independence theorem of 
Kinjo-Koseki \cite{KinjoKoseki} to identify $\widehat{\mathfrak{g}}^{C}(r,d)$ with $\widehat{\mathfrak{g}}^{C}(r,0)$, compatibly with 
the perverse filtration. This reduces the proof of the purity conjecture to the following equality, equivalent to \eqref{eq:intro1}:
\begin{equation}\label{eq:intro2}
	\sum_\ell \dim \gr^{\P}_{N-\ell}\; \widehat{\mathfrak{g}}^{C}_{2\ell}(r,0)q^\ell=A_{S_g,r}(q) .
\end{equation}
The advantage of replacing $\mathfrak{M}^{\ss}(r,d)$ with $\mathfrak{M}^{\ss}(r,0)$ is that, by Davison's neighborhood theorem, 
there is a local analytic isomorphism between a neighborhood of the trivial Higgs bundle 
$(\mathcal{O}_C^r,\phi=0)$ in $\mathfrak{M}^{\ss}(r,0)$ and a neighborhood of $0$ in the stack $\mathfrak{Rep}_r\Pi_{S_g}$ 
of dimension $r$ representations of the preprojective algebra of $S_g$. 
This is the bridge between the Higgs moduli space and the quiver variety. 
Note that the same bridge exists for the untwisted character variety,
however  the $\chi$-independence phenomenon has not been established on the Betti side. 
The remainder of the proof consists of the following three steps:

\medskip

\noindent 
\textit{Step $1$.} 
A recent general stacky form of the P=C theorem due to R. Picovschi \cite{Raphael}, valid for arbitrary pairs $(r,d)$,
implies that $\bigoplus_\ell \gr^{\P}_{N-\ell}\; \widehat{\mathfrak{g}}^{C}_{2\ell}(r,0)$ coincides with the graded of the subspace of 
$\widehat{\mathfrak{g}}^{C}(r,0)$ generated by the action of a family of tautological classes $\overline{c}_i([\text{pt}])$, $i \geqslant 2$, on the 
reduced virtual fundamental class $[\mathfrak{M}^{\ss,\rd}(r,0)]$. This subspace restricts, under the local isomorphism 
$$(\mathfrak{M}^{\ss}(r,0), \mathcal{O}_C^r) \cong (\mathfrak{Rep}_r\Pi_{S_g},0)$$ to the subspace generated by the action of the 
tautological classes $c_i(V)$, $i \geqslant 2$, on the reduced virtual fundamental class $[\mathfrak{Rep}^\rd_r\Pi_{S_g}]$. 
Note that, the stack  $\mathfrak{Rep}_r\Pi_{S_g}$ being conical, its Borel-Moore homology is the same as that of any analytic 
neighborhood of its apex $0$.
	
\medskip

\noindent 
\textit{Step $2$.} We use the following quiver analog of the Markman theorem for the BPS Lie algebras of symplectic surfaces proven in \cite{DHKSV}, which holds for any finite quiver. The same theorem is proved in \cite{DH} using an independent method.

\begin{theorem}\label{thm: Markman1}
Let $Q=(Q_0,Q_1)$ be any quiver and $d$ a dimension vector for $Q$. The subspace of tautological classes
$$\H^{\taut}_*(\mathfrak{Rep}_d\Pi_Q,\Q)=\mathbb{Q}[\ch_i(V_j)\,;\,i \geqslant 1, j \in Q_0] \cap [\mathfrak{Rep}^\rd_d\Pi_Q] $$
in $\H_*(\mathfrak{Rep}_d\Pi_Q,\Q)$ coincides with the BPS Lie algebra $\widehat{\mathfrak{g}}^{\Pi_Q}(d)$. 		
\end{theorem}
By Theorem~\ref{thm:Markman2}, the restriction to the neighborhood of $\mathcal{O}_C^r$ establishes a graded, surjective 
morphism 
$$\bigoplus_\ell \gr^{\P}_{N-\ell}\; \widehat{\mathfrak{g}}^{C}_{2\ell}(r,0)\twoheadrightarrow 
\widehat{\mathfrak{g}}^{\Pi_{S_g}}(r) /c_1(V) \cdot \widehat{\mathfrak{g}}^{\Pi_{S_g}}(r).$$ 
The graded dimension of the right hand side is equal to the right hand side of \eqref{eq:intro2}. 
This establishes the $\geqslant$ form of \eqref{eq:intro2}.
	
\medskip

\noindent 
\textit{Step $3$.}  We use the  Hard Lefschetz theorem to identify, in the coprime case,  
$$\bigoplus_\ell\gr^{\P}_{N-\ell}\, \H_{2N+2\ell}(M^D(r,d),\Q)\cong\bigoplus_\ell \gr^{\P}_{N+\ell}\,\H_{2N}(M^{D}(r,d),\Q).$$
Although we do not know how to compute the graded dimension of the latter, 
its {ungraded} dimension is equal to the constant term of the Kac polynomial $A_{g,r}(\underline{\sigma})$ counting absolutely indecomposable vector bundles of rank $r$ and degree $d$ by \cite[Thm.~1.1, Cor.~1.7]{S}.
Further, we have
	$$A_{g,r}(0)=A_{Sg,r}(1).$$ 
This provides the desired ungraded $\leqslant$ form of \eqref{eq:intro2} and allows us to conclude that \eqref{eq:intro2} holds. The proof of Theorem~\ref{thm:purityconj} follows.
\medskip

\medskip

\medskip\subsection{Some corollaries} Several consequences are worth pointing out: 

\begin{corollary}
	For any pair $(r,d)$ there exists a canonical isomorphism
	$$\bigoplus_\ell \gr^{\P}_{N-\ell}\; {\mathfrak{g}}^{C}_{2\ell}(r,d) \cong {\mathfrak{g}}^{\Pi_{S_g}}(r) .$$
\end{corollary}

\begin{proof}
This follows from the $\chi$-independence of $\bigoplus_\ell \P_{N-\ell}\; \widehat{\mathfrak{g}}^{C}_{2\ell}(r,d)$ together with the relations 
$$\gr^{\P}_{N-\ell}\; {\mathfrak{g}}^{C}_{2\ell}(r,d) = \gr^{\P}_{N-\ell}\; \widehat{\mathfrak{g}}^{C}_{2\ell}(r,d),$$ 
see \eqref{eq:isomgradedpureperverse},
and 
$$\widehat{\mathfrak{g}}^{\Pi_{S_g}}(r)/c_1(V)\cdot \widehat{\mathfrak{g}}^{\Pi_{S_g}}(r) \cong {\mathfrak{g}}^{\Pi_{S_g}}(r).$$ 
Note that, contrary to the case of a general symplectic surface, the $\chi$-independence isomorphism for $T^*C$ is canonical as $H^2(T^*C,\Q)$ is one-dimensional.
\end{proof}

\begin{corollary}
There exists a surjective Lie algebra homomorphism 
$$\widehat{\g}^{C,(0)}=\bigoplus_r \widehat{\g}^C(r,0) \to \bigoplus_r \widehat{\g}^{\Pi_{S_g}}(r)=\widehat{\g}^{\Pi_{S_g}} .$$
\end{corollary}
	 
\begin{corollary}
For coprime $(r,d)$, the perverse filtration on $\H_{2N}(M^D(r,d),\Q)$ is given by
$$\sum_\ell \dim \gr^{\P}_{N+\ell}\;\H_{2N}(M^{D}(r,d),\Q)q^\ell= A_{S_g,r}(q) .$$
\end{corollary}

Note that the $\chi$-independence allows us to extend this last corollary to arbitrary pair $(r,d)$, by replacing $\H_{2N}(M^{D}(r,d),\Q)$ by $\mathfrak{g}^C_{2N}(r,d)$, where $\mathfrak{g}^C(r,d)$ stands for the (non-affinized) BPS Lie algebra.

\medskip

\medskip\subsection{Parabolic case and the Hausel-Letellier-Rodriguez-Villegas purity conjecture} The purity conjecture~\ref{conj:HRV} has an 
extension to the case of a curve with punctures and prescribed generic monodromies around the punctures, see \cite{HLRV11}. 
This parabolic purity conjecture involves a comet-shaped quiver.
 The statement is very similar to the unpunctured case.
We expect that our proof may be adapted to this setting: both the P=W and curious Hard Lefschetz theorems are proved in the parabolic 
setting, see \cite{HMMS}.
The Poincar\'e polynomial of the parabolic character variety is computed, see \cite{MAnnals}.
What remains to be extended to the parabolic case are the $\chi$-independence and P=C statements. We plan to return to this in the future.

\medskip

\medskip\subsection{Plan of the paper.} Section 2 serves as a reminder about COHAs and BPS Lie algebras associated to quivers. The quiver analog of Markman's theorem is proved in \S3. In \S\S4,5 we introduce the necessary notions regarding the COHA and BPS Lie 
algebra associated to Higgs moduli, and recall the main results of \cite{S17}, \cite{KinjoKoseki}, \cite{DHKSV} and \cite{Raphael}. After 
recalling Davison's neighborhood theorem, the associated restriction map is defined and studied in \S6. Finally, the proof of 
Theorem~\ref{thm:purityconj} is given in \S7. The appendix provides some detail about the COHA coproduct and the proofs of various 
results around the primitivity of tautological classes, adapted from \cite{DHKSV}.

\medskip

%=====================================================================
\medskip\section{Quivers, COHAs and the BPS Lie algebra of a quiver.}
%=====================================================================

\medskip

In this section, we recall the main properties of cohomological Hall algebras associated to quivers and their BPS Lie algebras, postponing 
some proofs and details to Appendix A. In addition to the usual COHA of the preprojective algebra $\Pi_Q$ of a quiver $Q$, we 
also consider the COHA of the endomorphism algebra $\C Q[\theta]$, see below. The two COHAs are related by dimensional reduction.
The latter is better suited to our proof of Markman's tautological generation theorem for the BPS Lie algebra, 
see Theorem~\ref{thm:Markman2}. Throughout the entire paper, we work over $\C$.

\medskip

\medskip\subsection{Stacks of quiver representations} Let $Q=(Q_0,Q_1)$ be an arbitrary finite quiver, with possible multiple arrows and edge loops. The Euler form $\langle-,-\rangle$ and the symmetrized Euler form $(-,-)$ are given by
\begin{equation}\label{eq:Eulerform}
\langle d,e\rangle=\sum_i d_ie_i -\sum_{a\colon i\to j}d_ie_j, \quad (d,e)=2\sum_i d_ie_i-\!\sum_{a\colon i\to j}(d_ie_j+ d_je_i).
\end{equation}
Let $\overline{Q}=(Q_0, Q_1 \cup Q_1^*)$ be the double quiver. For a dimension vector $d\in\N Q_0$ we put
$$E_d=\bigoplus_{a \in Q_1} \text{Hom}(\C^{d_{s(a)}}, \C^{d_{t(a)}}), \quad \overline{E}_d=T^*E_d, \quad
GL_d=\prod_{i\in Q_0}GL_{d_i},\quad \gl_d=\operatorname{Lie}(GL_d) .$$
For $d=(d_i)_{i \in Q_0}$ we set
$$|d|=\sum_i d_i, \quad \tr: \gl_d \to \C,\; \tr((\theta_i)_i)=\sum_i \tr(\theta_i) .$$
The stacks of $d$-dimensional representations of $Q$ and $\overline{Q}$ are
$$\Rep_dQ=E_d/GL_d, \quad \Rep_d\overline{Q}=\overline{E}_d/GL_d .$$
The stack $\Rep_dQ$ is smooth of dimension $-\langle d,d\rangle$.
The preprojective algebra $\Pi_Q$ of $Q$ is the quotient of $\C \overline{Q}$ by the two-sided ideal generated by $\sum_{a \in Q_1} [x_a, x_a^*]$. The moment map for the Hamiltonian action of $GL_d$ on $T^*E_d$ is 
\begin{equation}\label{eq:momentmap}
\mu_d\colon T^*E_d\to\gl_d,\quad (x,x^*)\mapsto\sum_{a\in Q_1}[x_a,x^*_a].
\end{equation}
The derived stack of $d$-dimensional $\Pi_Q$-modules is thus
$$
\Rep_d\Pi_Q=T^*\big(\Rep_dQ)=\mu_d^{-1}(0)/GL_d,
$$
where $\mu_d^{-1}(0)$ is the derived zero locus. It is $0$-shifted symplectic, quasi-smooth, of virtual dimension
$$
\vdim\Rep_d\Pi_Q=2\dim\Rep_dQ=-(d,d).
$$ 
We  also consider the additive version of the inertia stack of $\Rep Q$, i.e.,  the stack of $\C Q[\theta]$-representations. 
It is the quotient by $GL_d$ of a derived zero locus in an affine space
$$\Rep_d\IQ= T[-1]\Rep_dQ=\{(x,\theta)\;;\;x\in\Rep_dQ,\ \theta\in\gl_d,\ [\theta,x]=0\}/GL_d$$
It is quasi-smooth of virtual dimension $0$. 
Finally, the critical locus of the tripled quiver $\widetilde{Q}$ 
with potential 
\begin{align}\label{potential}W=\sum_{a\in Q_1}x^*[\theta,x]\end{align}
is the stack of modules over the algebra $\Pi_Q[\theta]$
$$\tRep_d\Pi_Q=T^*[-1]\Rep_d\Pi_Q=\{(x,x^*,\theta)\;;\;\mu_d(x,x^*)=0,\ [\theta,x]=[\theta,x^*]=0\}/GL_d.$$
It is a $(-1)$-shifted symplectic stack of virtual dimension $0$, carrying a canonical orientation,
see \cite[ex.~3.20]{KPS}. Forgetting $\theta$ or $x^*$ gives the two dimensional reductions morphisms
$$
\begin{tikzcd}[column sep=small]
& \tRep_d\Pi_Q \arrow[dl,"\kappa'"'] \arrow[dr,"\kappa"] & \\
\Rep_d\IQ & & \Rep_d\Pi_Q
\end{tikzcd}
$$
All three stacks map to $BGL_d$ compatibly.

\medskip

\medskip\subsection{Conical substacks}
Inside $\tRep_d\Pi_Q$ and $\Rep_d\IQ$ we consider the closed substacks defined by the following conditions on $\theta$
$$\begin{array}{lll}
\Lambda_d^0,&\ \Rep^0_d\IQ& :\ \theta=0,\\[1mm]
\Lambda_d^\nil,&\ \Rep^\nil_d\IQ &:\ \theta\ \text{nilpotent},\\[1mm]
\Lambda_d^\scal,&\ \Rep^\scal_d\IQ& :\ \exists\, c\in\C,\ \theta=c\cdot\mathrm{id},\\[1mm]
\Lambda_d^\snil,&\ \Rep^\snil_d\IQ &:\ \exists\, c\in\C,\ \theta-c\ \text{nilpotent}.
\end{array}
$$
For each $\star \in \{0,\nil,\scal,\snil\}$ we have
$$\Lambda_d^\star=(\kappa')^{-1}(\Rep^\star_d\IQ).$$ 
There are closed immersions 
$$\begin{tikzcd}[column sep=small]
	 \Lambda_d^0 \arrow[d] \arrow[r] & \Lambda_d^\nil \arrow[d] \\
	\Lambda_d^\scal \arrow[r]& \Lambda_d^\snil
\end{tikzcd}$$ 
and similarly for $\Rep^\star_d\IQ$.
The stacks $\Lambda_d^\star$ and $\Rep^\star_d\IQ$ are taken with their {classical}, reduced structure. Observe that $\Rep^0_d\IQ=\Rep_dQ$. The following lemma is obvious.

\begin{lemma}\label{lem:translation}
For $d\neq0$, the map 
$$(x,\theta)\mapsto\big((x,\theta-c\,\mathrm{id}),c\big),\quad c=\tr(\theta)/|d|$$ 
gives $GL_d$-equivariant isomorphisms
\[
\Rep^\snil_d\IQ\cong\Rep^\nil_d\IQ\times\A^1,\quad\Rep^\scal_d\IQ\cong\Rep^0_d\IQ\times\A^1 .
\]
They are compatible with the inclusions and with the maps to $BGL_d$.
\qed
\end{lemma}

\medskip

\medskip\subsection{COHAs of quivers}
 For $l \in \Z$ and $d \in \N Q_0$ we put
 \begin{gather*}
 \bfA_{\Pi_Q,l}(d)=\H_{l-(d,d)}(\Rep_d\Pi_Q,\Q),\quad \bfA_{\Pi_Q}(d)=\bigoplus_l \bfA_{\Pi_Q,l}(d), 
 \quad \bfA_{\Pi_Q}=\bigoplus_d \bfA_{\Pi_Q}(d),\\
 \bfA_{IQ,l}(d)=\H_{l}(\Rep_d IQ,\Q), \quad \bfA_{IQ}(d)=\bigoplus_l \bfA_{IQ,l}(d), \quad \bfA_{IQ}=\bigoplus_d \bfA_{IQ}(d)
 \end{gather*}
 and likewise for $\bfA^{\nil}_{IQ}$. 
 All three of $\bfA_{\Pi_Q}, \bfA_{IQ}$ and $\bfA_{IQ}^{nil}$ are equipped with $\N Q_0 \times \Z$-graded 
 cohomological Hall algebra (COHA) structures. 
 The multiplication in $\bfA_{\Pi_Q}$ comes from the correspondence
 \begin{equation}\label{diag:defprodcoha1}
 \begin{tikzcd}[column sep=small]
 	\Rep_d\Pi_Q \times \Rep_{d'}\Pi_Q && \Filt_{d,d'}\Pi_Q \arrow[ll,"\gr_{d,d'}"'] \arrow[rr,"\ev_{d,d'}"] && \Rep_{d+d'}\Pi_Q ,
 \end{tikzcd}
\end{equation}
 where the stack $\Filt_{d,d'}\Pi_Q$ classifies short exact sequences 
 $$0 \to M_{d'} \to M \to M_d \to 0$$ 
 of $\Pi_Q$-modules of respective dimensions $d,d'$ and $d+d'$, and where the maps $\gr_{d,d'}$ and $\ev_{d,d'}$ associate to such a sequence its extreme, resp. middle terms. The map $\gr_{d,d'}$ is quasi-smooth while $\ev_{d,d'}$ is proper. This allows one to define an associative product 
 $$\bfA_{\Pi_Q}(d) \otimes \bfA_{\Pi_Q}(d') \to \bfA_{\Pi_Q}(d+d'), \quad c \otimes c' \mapsto (\ev_{d,d'})_*(\gr_{d,d'})^!(c \boxtimes c')$$
 where $\gr^!_{d,d'}$ is the virtual pullback. With the shift by $-(d,d)$ in the definition of $\bfA_{\Pi_Q,l}(d)$, this multiplication is of cohomological degree zero. The construction of the algebra structure on $\bfA_{IQ}$ and $\bfA_{IQ}^{\nil}$ is entirely similar.  
 Moreover, denoting by 
 \begin{align}\label{iota-nil}\iota^{\nil}:  \Rep^{\nil}IQ\to\Rep IQ
 ,\quad
 \Rep^{\nil}IQ=\bigcup_d\Rep^{\nil}_dIQ
 ,\quad
 \Rep IQ=\bigcup_d\Rep_dIQ\end{align}
 the closed embedding, the pushforward 
 $$\iota^\nil_*: \bfA^\nil_{IQ} \to \bfA_{IQ}$$ is an algebra homomorphism. 
 We refer to \cite{SVIHES}, \cite{KV19} and \cite{PortaSala} for details. We also define
 $$\bfA_{IQ,l}^{\snil}(d)=\H_l(\Rep^\snil_dIQ,\Q), \quad \bfA_{IQ,l}^{\scal}(d)=\H_l(\Rep^\scal_dIQ,\Q), \quad \bfA_{IQ,l}^{0}(d)=\H_l(\Rep^0_dIQ,\Q)=\H_l(\Rep_dQ, \Q)$$
 Note that 
 $$\bfA^{\snil}_{IQ}=\bigoplus_{l,d} \bfA_{IQ,l}^{\snil}(d)
 ,\quad
 \bfA^{\scal}_{IQ}=\bigoplus_{l,d} \bfA_{IQ,l}^{\scal}(d)
 ,\quad
 \bfA^{0}_{IQ}=\bigoplus_{l,d} \bfA_{IQ,l}^{0}(d)$$ do not carry natural algebra structures as they 
 correspond to substacks that are not stable under extensions. 
 There is a $1d$-COHA structure on $\bfA^0_{IQ}$. We do not consider it here.

 \medskip
 
 Both $\bfA_{\Pi_Q}$ and $\bfA_{IQ}$ come from the same critical COHA associated to $\tRep \Pi_Q$. 
 By \cite{BBDJS15}, the stack $\tRep_d\Pi_Q$ carries a Joyce perverse sheaf $\varphi_{\tRep_d\Pi_Q}$ and
 $$\widetilde{\bfA}=\bigoplus_{l,d} \H^{-l}\left(\tRep_d\Pi_Q,\varphi_{\tRep_d\Pi_Q}\right)$$
 carries an associative algebra structure defined via an induction diagram very similar to \eqref{diag:defprodcoha1}. Moreover, there are dimensional reduction isomorphisms
$$\kappa_*\varphi_{\tRep_d\Pi_Q}\cong \DQ_{\Rep_d\Pi_Q}\otimes\L^{-\frac{1}{2}(d,d)} \quad \text{and} \quad 
%\DQ_{\Rep_d\Pi_Q}\otimes\L^{-\vdim \Rep_d\Pi_Q}=
\kappa'_*\varphi_{\tRep_d\Pi_Q}\cong
%\DQ_{\Rep\IQ}\otimes\L^{-\vdim \Rep_dIQ}=
\DQ_{\Rep\IQ}$$ 
see \cite{D17} or  \cite[\S2]{DHKSV}. Taking global sections gives the isomorphisms
\begin{equation}\label{eq:star}
	\bfA_{\Pi_Q,l}(d)\cong \widetilde{\bfA}_l(d) \cong \bfA_{IQ,l}(d) .
\end{equation}
These isomorphisms are $\H_{GL_d}^*$-linear and intertwine the three COHA products.
 %{{\color{red}{In addition, they are equivariant for the torus $\T=\GM$ acting with weight $1$ on $\theta$, weight $-1$ on $x^*$ and weight $0$ on $x$.This action preserves the potential $W$. [Is this used anywhere ??]}}}
For the closed conical substacks $\Rep^\star\IQ$ with $\star \in \{\nil,\snil,\scal,0\}$,
a base change with respect to the map $\kappa'$ yields in the same fashion isomorphisms
\begin{equation}\label{eq:star-supp}
	\H^{\Lambda^\star}_*(\Rep_d\Pi_Q,\Q)\cong\H_*(\Rep^\star_d\IQ,\Q),
\end{equation}
which are compatible with \eqref{eq:star} via the pushforward along the closed immersions 
$$\Rep^\star_d\IQ\hookrightarrow\Rep_d\IQ.$$
The left hand side of \eqref{eq:star-supp} is the microlocal homology with support.
It is defined by
$$ \H^{\Lambda^\star}_l(\Rep_d\Pi_Q,\Q)=
\H^{-l}\big(\tRep_d\Pi_Q\,,\, \tilde{\iota}^\star_! (\tilde{\iota}^\star)^!\varphi_{\tRep_d\Pi_Q}\big) ,$$
where $\tilde{\iota}^\star: \Lambda^\star_d \to \tRep_d\Pi_Q$ is the closed immersion, see \cite[\S 2.3]{DHKSV}. 
In particular, there is an algebra isomorphism
$$\bigoplus_d \H^{\Lambda^{\nil}}(\Rep_d\Pi_Q\,,\,\Q) \cong \bfA^\nil_{IQ} .$$

\medskip\subsection{Purity} For future use, we mention the following result.
\begin{proposition}\label{prop: purity}
The mixed Hodge structures on the spaces
$$\bfA_{\Pi_Q}(d),\; \bfA_{IQ}(d),\; \bfA^{\star}_{IQ}(d), \quad \star \in \{0,\scal,\nil,\snil\}$$
are all pure, for all $d \in \N Q_0$.	
\end{proposition}
\begin{proof}
Purity of $\bfA_{\Pi_Q}(d)$ is shown in \cite{DavisonKac}, \cite{D24}. The same thus holds for $\bfA_{IQ}$. Purity of the conical substacks $\Rep_d^\star IQ$ is much easier: it is obvious for $\star\in\{0,\scal\}$ and the case of $\star=\snil$ follows from the case $\star=\nil$. In the latter case, it is a consequence of the Jordan stratification, see, e.g., \cite{MS20} or Corollary~\ref{cor:purity-nil} below. 
\end{proof}

\medskip

\medskip\subsection{The BPS sheaves and the BPS Lie algebras $\widehat{\g}^{\Pi_Q}$}
We next briefly recall the main features of the theory of BPS sheaves and Lie algebras, as developed in the work of Davison and his
collaborators in \cite{DM}, \cite{DHSM2}, \cite{DHKSV}. We refer to Appendix~A for some details. Denote by 
$$p: (\Rep_d\Pi_Q)^\red \to Rep_d\Pi_Q=\mu_d^{-1}(0)/\!\!/GL_d $$
the good moduli spaces and by 
$$\oplus_{d,d'}: Rep_d\Pi_Q \times Rep_{d'}\Pi_Q \to Rep_{d+d'}\Pi_Q$$
the addition maps. There are pure mixed Hodge modules $\mathcal{BPS}^Q_d$ over $Rep_d\Pi_Q$ such that
\begin{equation}\label{eq:integrality}
\text{Sym}_{\oplus}\Big( \bigoplus_d \mathcal{BPS}^{\Pi_Q}(d) \otimes\\H_{\GM}^*\Big) \cong \bigoplus_
d p_*\left(\DQ_{\Rep_d\Pi_Q} \otimes \mathbb{L}^{-\frac{1}{2}(d,d)}\right)
\end{equation}
 The BPS Lie algebras and affinized BPS Lie algebras of the quiver $Q$ are defined as
\begin{align*}
\g^{\Pi_Q}&=\bigoplus_d \g^{\Pi_Q}(d)= \bigoplus_d \H^{-*}\Big(Rep_d\Pi_Q, \mathcal{BPS}^{\Pi_Q}(d)\Big),\\
\widehat{\g}^{\Pi_Q}&=\bigoplus_d \widehat{\g}^{\Pi_Q}(d)= \bigoplus_d \H^{-*}\left(Rep_d\Pi_Q, \mathcal{BPS}^{\Pi_Q}(d)\otimes H_{\GM}^*\right).
\end{align*}
They are $\N Q_0\times\Z$-graded subspaces in $\bfA_{\Pi_Q}$.
The generator $u \in \H_{\GM}^2$ is in degree $-2$.
Both are Lie subalgebras of the COHA. Moreover, we have
\begin{equation}\label{eq:A=Ug}
 \text{U}\left(\widehat{\g}^{\Pi_Q}\right) = \bfA_{\Pi_Q}
\end{equation}
Passing to global sections in \eqref{eq:integrality} recovers the PBW isomorphism
\begin{equation}\label{eq:PBW1}
 \Sym \Big( \bigoplus_d \widehat{\g}^{\Pi_Q}(d)\Big) \xrightarrow{\cong} \bfA_{\Pi_Q}
\end{equation}
where the map is given by the COHA multiplication.
Finally, there is a comultiplication \cite{DHKSV}
\begin{equation}\label{eq:coproduct}
\Delta: \bfA_{\Pi_Q} \to \bfA_{\Pi_Q} \otimes \bfA_{\Pi_Q}
\end{equation}
which is compatible with the COHA product and whose space of primitive elements is $\widehat{\g}^{\Pi_Q}$. 
By \cite{DavisonKac} the graded dimensions of $\g^{\Pi_Q}(d)$ and $\widehat{\g}^{\Pi_Q}(d)$ are given by
\begin{equation}\label{eq:dimBPSLiealg}
\sum_{l} \dim \g^{\Pi_Q}_l(d)t^l= A_{Q,d}(t^2), \quad 
\sum_{l} \dim \widehat{\g}^{\Pi_Q}_l(d)t^l=\frac{1}{1-t^{-2}} A_{Q,d}(t^2) ,
\end{equation}
where $A_{Q,d}(q)$ is Kac's polynomial counting absolutely indecomposable representations of $Q$ of dimension $d$ 
over the finite field $\mathbb{F}_q$. 
We define the Lie subalgebras 
$$\g^{IQ}\,,\,\widehat{\g}^{IQ}\subset\bfA_{IQ}$$ 
through the identification \eqref{eq:star}, so that \eqref{eq:A=Ug}, \eqref{eq:PBW1} hold with $IQ$ in place of $Q$. 

\begin{remark}
The above BPS Lie algebra formalism also applies to  
$$\H^{\Lambda^{\nil}}(\Rep\Pi_Q,\Q) \cong \bfA^\nil_{IQ}$$ 
There is a Lie subalgebra 
$$\widehat{\g}^{IQ,\nil} \subset \bfA_{IQ}^\nil$$ such that 
$$\text{U}(\widehat{\g}^{IQ,\nil})=\bfA_{IQ}^\nil,\quad
\widehat{\g}^{IQ,\nil} \cong \widehat{\g}^{IQ}[-2] \text{ (as vector spaces)}$$
where $[n]$ indicates a 
shift by $n$ in the cohomological grading. The Lie algebra structures on $\widehat{\g}^{IQ,\nil}$ and  $\widehat{\g}^{IQ}$ are very 
different: the latter is the positive half of a generalized Kac-Moody algebra by \cite{DHSM2} (see also \cite{DHSM1}), while the former can be shown to be abelian by \cite{NS}. We will 
not need this.   
\end{remark}

\medskip

The following characterizations of $\widehat{\g}^{\Pi_Q}$ and $\widehat{\g}^{IQ}$ play a key role in this paper. Recall that
\begin{align}\label{iota-snil}\tilde{\iota}^{\snil}: \Lambda^{\snil}_d \to \tRep_d\Pi_Q, \quad  \iota^\snil\colon\Rep^\snil_d\IQ\to\Rep_d\IQ
\end{align}
denote the closed immersions. 

\begin{proposition}\label{prop:seminilpcharacterizationBPS}
Under the dimensional reduction isomorphisms \eqref{eq:star} and \eqref{eq:star-supp} we have
\begin{equation}\label{eq:supportcondition1}
	\Im\Big(\tilde{\iota}^\snil_*\colon\H^{\Lambda^\snil}(\Rep_d\Pi_Q\,,\,\Q) \to\bfA_{\Pi_Q}(d)\Big)=\widehat{\g}^{\Pi_Q}(d) ,
\end{equation}
\begin{equation}\label{eq:supportcondition2}
\Im\Big(\iota^\snil_*\colon\bfA_{IQ}^\snil(d)\to\bfA_{IQ}(d)\Big)=\widehat{\g}^{IQ}(d) .
\end{equation}
\end{proposition}
\begin{proof}Equality \eqref{eq:supportcondition2} follows from \eqref{eq:supportcondition1} by base change and dimensional reduction. The proof of \eqref{eq:supportcondition1} is given in \cite{DHKSV} in a more general setting. We provide some details in Appendix~A.
\end{proof}

\medskip

\medskip\subsection{Tautological classes} 
For each vertex $i\in Q_0$ let $V_i$ be the tautological rank $d_i$ bundle on $\Rep_d\Pi_Q$. We will use the same notation for $\Rep_dIQ$. 

\begin{lemma}\label{lem:all-tautological}
	\hfill
	\begin{enumerate}[label=$\mathrm{(\alph*)}$,leftmargin=8mm,itemsep=1mm]
		\item
		The pullbacks $\H_{GL_d}^*\to\H^*(\Rep_d\Pi_Q,\Q)$ and $\H_{GL_d}^*\to\H^*(\Rep_dIQ,\Q)$ are isomorphisms.
		\item
		$\H^*(\Rep_d\Pi_Q,\Q)$ and $\H^*(\Rep_dIQ,\Q)$ are generated by the Chern classes of the $V_i$.
	\end{enumerate}
\end{lemma}

\begin{proof}
	The cohomology of a derived stack being equal to that of its classical truncation,
	$$\H^*(\Rep_d\Pi_Q,\Q)=\H^*_{GL_d}(\mu_d^{-1}(0)^\cl,\Q).$$ 
	The scaling action $t\cdot(x,x^*)=(tx,tx^*)$ on $\Rep_d\overline{Q}$ commutes with $GL_d$ and preserves the conical subvariety 
	$\mu_d^{-1}(0)^\cl$.
	Hence $\mu_d^{-1}(0)^\cl$ is $GL_d$-equivariantly contractible to $0$, and we get
	\begin{align}\label{taut0}\H^*_{GL_d}(\mu_d^{-1}(0)^\cl,\Q)\cong\H_{GL_d}^*.\end{align}
	The proof is the same for $\Rep_dIQ$. The second statement follows obviously from the first.
\end{proof}

Consider the graded ring 
$$S_Q=\Q[\ch_{i,j}\;;\; j \in Q_0\,,\,i \geqslant 1],\quad \deg(\ch_{i,j})=2i\ .$$
For any $d\in \N Q_0$ the assignement $\ch_{i,j} \mapsto \ch_i(V_j)$ defines surjective algebra morphisms
$$a_d: S_Q \to \H^*(\Rep_d\Pi_Q,\Q), \quad a_d: S_Q \to \H^*(\Rep_dIQ,\Q)\ .$$
The action of cohomology on Borel-Moore homology gives rise to actions of $S_Q$ on $\bfA_{\Pi_Q}(d)$ and $\bfA_{IQ}(d)$ by 
$$c\bullet x=a_d(c)\cap x.$$ 
To unburden the notation, we will often suppress the notation $a_d$.

\begin{lemma}\label{lem:taut_classes_intertwine}
	The isomorphism \eqref{eq:star} intertwines the actions of $S_Q$ on $\bfA_{\Pi_Q}$ and $\bfA_{IQ}$.
	\qed
\end{lemma}
 
The main compatibility properties between the action of $S_Q$ and the COHA operations are summarized in the following lemma. 
We equip $S_Q$ with the unique comultiplication $\Delta$ satisfying
$$\Delta(\ch_{i,j})=\ch_{i,j} \otimes 1 + 1 \otimes  \ch_{i,j} ,\quad \forall\,j \in Q_0\,,\, i \geqslant 1 .$$

\begin{lemma}\label{lem:hopf-module}
	For $c\in S_Q$ and homogeneous $x,y\in\bfA_{\Pi_Q}$, resp. $x,y \in \bfA_{IQ}$, we have
	\begin{enumerate}[label=$\mathrm{(\alph*)}$,leftmargin=8mm,itemsep=1mm]
		\item $c\bullet(xy)=\sum(c_{(1)}\bullet x)(c_{(2)}\bullet y)$,
		\item $\Delta(c\bullet x)=\Delta(c)\bullet\Delta(x)$,
		\item $S_Q\bullet \widehat{\g}^{\Pi_Q}\subseteq\widehat{\g}^{\Pi_Q}$ and likewise for $\widehat{\g}^{IQ}$.
	\end{enumerate}
\end{lemma}
\begin{proof}
The proof is the same as in \cite[\S 6]{DHKSV} (see Appendix A. for a brief sketch of the argument).
\end{proof}

Although the Lie subalgebra $\g^{\Pi_Q}$ is defined canonically by geometric means, it may also be recovered from 
$\widehat{\g}^{\Pi_Q}$ as follows. Assuming $d\neq 0$, the tautological line bundle 
$$\mathcal{L}_d=\bigotimes_i \det V_i$$ 
has a nontrivial weight under the standard $\BGM$-action by scalar multiplication.
It is thus a global equivariant parameter in the sense of \cite[\S2.5.2]{DHKSV}. 
Then the inclusion $\g^{\Pi_Q} \subset \widehat{\g}^{\Pi_Q}$ induces an isomorphism
\begin{equation}\label{eq:gfromghat}
\g^{\Pi_Q} \cong \widehat{\g}^{\Pi_Q}\Big/ \Big( \bigoplus_d c_1(\mathcal{L}_d)\bullet \widehat{\g}^{\Pi_Q}(d)\Big).
\end{equation}
Using the relation $\mathcal{L}_d=\gr_{d',d''}^*(\mathcal{L}_{d'} \otimes \mathcal{L}_{d''})$ whenever $d=d'+d''$, one can show that 
$\bigoplus_d c_1(\mathcal{L}_d)\bullet \widehat{\g}^{\Pi_Q}(d)$ is a Lie ideal and \eqref{eq:gfromghat} is a Lie algebra isomorphism.

\medskip

\medskip\subsection{The reduced virtual
	 fundamental class} The moment map $\mu_d$ of \eqref{eq:momentmap} factors through 
$$\gl_d^0=\ker\big(\tr:\gl_d\rightarrow\C\big).$$ 
Let $\mu^0_{d}\colon{\Rep}_d\overline{Q}\to\gl_d^0$ be the factorization.
Set $$\Rep^\rd_d\Pi_Q=(\mu_d^0)^{-1}(0)/GL_d .$$
Then, we have 
$$(\Rep^\rd_d\Pi_Q)^\cl=(\Rep_d\Pi_Q)^{\cl}
,\quad
\vdim(\Rep^\rd_d\Pi_Q)=\vdim\Rep_d\Pi_Q+1=1-(d,d)$$
and the isomorphism of  derived stacks
$$\Rep_d\Pi_Q\cong\Rep^\rd_d\Pi_Q\times\Spec\Sym(\C[1])$$
%This is the reduction of \cite[\S2.3.4]{DHKSV} for the central $\BGM\subset \B G_d$.
%The moment map of the scalar $\BGM$-action is $\sum\tr\circ\mu$. 
The reduced virtual fundamental class of $\Rep_d\Pi_Q$ is 
$$[\Rep_d\Pi_Q]^\rd=[\Rep^\rd_d\Pi_Q]\in\H_{-2(d,d)+2}(\Rep_d\Pi_Q,\Q).$$
Note that $[\Rep_d\Pi_Q]=0$ for $d\neq0$, because the obstruction theory contains a summand $\C$.

We need the following microlocal characterization of $[\Rep^\rd_d\Pi_Q]$. 
Namely, by \cite[\S7.5]{DHKSV}, the reduced virtual class
$[\Rep^\rd_d\Pi_Q]$ is the image of $1 \in \Q$ under the chain of maps
$$
\Q\cong\H^{(d,d)-2}\big(\tRep_d\Pi_Q,\,(\tilde{\iota}^\scal)_!(\tilde{\iota}^\scal)^!\varphi_{\tRep_d\Pi_Q}\big)\to
\H^{(d,d)-2}\big(\tRep_d\Pi_Q,\,\varphi_{\tRep_d\Pi_Q}\big)\cong\H_{2-2(d,d)}(\Rep_d\Pi_Q,\Q).
$$

\begin{proposition}\label{prop:reducedvirtuafund_is_primitive}
	For any $d \in \N Q_0$ we have $[\Rep^\rd_d\Pi_Q] \in \widehat{\g}^{\Pi_Q}(d)$.
\end{proposition}

\begin{proof} This follows from the microlocal characterization of $\widehat{\g}^{\Pi_Q}$ and the chain of inclusions
	$$\Lambda^\scal_d \subset \Lambda^\snil_d \subset \tRep_d\Pi_Q .$$  
\end{proof}

\begin{proposition}\label{prop:red-class}
	Under the isomorphism \eqref{eq:star}, the reduced virtual fundamental class
	$[\Rep^\rd_d\Pi_Q]$ is the pushforward $(\iota^\scal)_*[\Rep^\scal_d\IQ]$
	of the fundamental class of the smooth stack 
	$\Rep^\scal_d\IQ\cong\Rep_dQ\times\A^1$. 
\end{proposition}

\begin{proof}
	The square
	$$
	\begin{tikzcd}
		\Lambda_d^\scal \arrow[r,hook,"\tilde{\iota}^\scal"] \arrow[d,"\kappa'"'] & \tRep_d\Pi_Q \arrow[d,"\kappa'"]\\
		\Rep^\scal_d\IQ \arrow[r,hook,"\iota^\scal"] & \Rep_d\IQ
	\end{tikzcd}
	$$
	is Cartesian. Since $\kappa'$ is proper on the relevant supports, base change gives a commutative diagram
	$$
	\begin{tikzcd}
		\kappa'_*(\tilde{\iota}^\scal)_!(\tilde{\iota}^\scal)^!\varphi_{\tRep_d\Pi_Q} \arrow[r] \arrow[d,"\cong"'] & \kappa'_*\varphi_{\tRep_d\Pi_Q}\arrow[d,"\cong"]\\
		(\iota^\scal)_!\,(\iota^\scal)^!\big(\DQ_{\Rep_d\IQ}\big)
		=(\iota^\scal)_!\DQ_{\Rep_d^\scal\IQ} \arrow[r] & \DQ_{\Rep_d\IQ},
	\end{tikzcd}
	$$
	in which the horizontal arrows are the counit maps. On global sections the bottom arrow is the map
	$$(\iota^\scal)_*\colon\H_*^\BM(\Rep^\scal_d\IQ,\Q)\to\H_*^\BM(\Rep_d\IQ,\Q).$$
	The top degree of $\H_*^\BM(\Rep^\scal_d\IQ,\Q)$ 
	is spanned by the fundamental class, which corresponds to the generator of the top left-hand side.
\end{proof}

Thanks to Lemma~\ref{lem:taut_classes_intertwine} and Proposition~\ref{prop:red-class}, the isomorphism \eqref{eq:star} identifies 
\[
c\cap[\Rep^\rd_d\Pi_Q]\ \xleftrightarrow{\eqref{eq:star}}\ (\iota^\scal)_*\big(c\cap[\Rep^\scal_d\IQ]\big)
,\quad \forall c\in S_Q.
\]

We define
$$\widehat{\g}^{\Pi_Q}_{\taut}(d)=S_Q \bullet [\Rep^\rd_d\Pi_Q], \quad \widehat{\g}^{\Pi_Q}_\taut=\bigoplus_d \widehat{\g}^{\Pi_Q}_{\taut}(d) $$
and likewise for $\widehat{\g}^{IQ}_\taut$. By Lemma~\ref{lem:hopf-module} (c), we have
\begin{equation}\label{eq:obviousinclusion}
\widehat{\g}^{\Pi_Q}_\taut \subseteq \widehat{\g}^{\Pi_Q}, \quad  \widehat{\g}^{IQ}_\taut \subseteq \widehat{\g}^{IQ} .
\end{equation}

%%%%%%%%%%%%%%%%%%%%%%%%%

\medskip

%=====================================================================
\medskip\section{A quiver analog of Markman's theorem.}
%=====================================================================

\medskip

The aim of this section is to prove the following quiver analog of \cite[Thm.~7.5]{DHKSV}.

\begin{theorem}\label{thm:Markman2}
We have
\begin{equation}\label{eq:Markmanequality}
	\widehat{\g}^{\Pi_Q}_\taut = \widehat{\g}^{\Pi_Q}, \quad  \widehat{\g}^{IQ}_\taut = \widehat{\g}^{IQ} ,
\end{equation}
Equivalently, the inclusions of \eqref{eq:obviousinclusion} are isomorphisms.
\end{theorem}

The two equalities in \eqref{eq:Markmanequality} are equivalent to each other. 
By the PBW isomorphism \eqref{eq:PBW1}, it is enough to prove that the 
COHA multiplication induces either of the equivalent isomorphisms
\begin{equation}\label{eq:reduction1}
\Sym\Big(\bigoplus_d\; \widehat{\g}^{\Pi_Q}_{\taut}(d)\Big)\xrightarrow{\cong}\bigoplus_d\bfA^{\Pi_Q}(d)
, \quad 
\Sym\Big(\bigoplus_d\; \widehat{\g}^{IQ}_{\taut}(d)\Big)\xrightarrow{\cong}\bigoplus_d\bfA^{IQ}(d) .
\end{equation}

Our strategy is to prove the second isomorphism in \eqref{eq:reduction1}. To this end, we consider the nilpotent COHA $\bfA^\nil_{IQ}$, the 
seminilpotent vector space $\bfA^\snil_{IQ}$ and we use the microlocal characterisation of $\widehat{\g}^{IQ}$ in 
Proposition~\ref{prop:seminilpcharacterizationBPS}.

\begin{remark} Theorem~\ref{thm:Markman2} should be compared to the Kirwan surjectivity theorem for Nakajima quiver varieties which is proved in \cite{MN}. Here, we are considering the entire stack of representations of the preprojective algebra and not an open substack of semistable points for some stability condition, which accounts for the fact the entire $\H_*(\Rep_d\Pi_Q,\Q)$ is not tautologically generated.
\end{remark}

\medskip

\medskip\subsection{The map $\Phi$.}
Let us put
\begin{align*}
	\widehat{\g}^0(d)&=S_Q\bullet[\Rep^0_d\IQ]\subseteq\bfA_{IQ}^\nil(d) .
\end{align*}
For $d \neq 0$ we consider the composition of maps
$$\Phi\colon\bfA_{IQ}^\nil(d)\xrightarrow{\cong} \bfA^{\snil}_{IQ}(d) \xrightarrow{\iota^{\snil}_*} \bfA_{IQ}(d)$$ 
The first map is the smooth pullback along the projection of Lemma~\ref{lem:translation}. 
We also set $\Phi(1)=1$ in degree $0$. Note that $\Phi$ is of (co)homological degree $2$.

\begin{proposition}\label{prop:Phi}
	The following hold:
	\begin{enumerate}[label=$\mathrm{(\alph*)}$,leftmargin=8mm,itemsep=1mm]
		\item $\Phi(\bfA_{IQ}^\nil(d))=\widehat{\g}^{IQ}(d)$.
		\item $\Phi(\widehat{\g}^0(d))=\widehat{\g}^{IQ}_{\taut}(d)$.
		\item $\Phi(x*x')=0$ for all $x\in\bfA^\nil_{IQ}(d)$, $x'\in\bfA^\nil_{IQ}(d')$ with $d,d'\neq0$, where $*$ is the COHA product on $\bfA^\nil_{IQ}$.
	\end{enumerate}
\end{proposition}

\begin{proof} The first statement is a direct consequence of Proposition~\ref{prop:seminilpcharacterizationBPS} and the fact that 
$$\Rep_d^\snilp IQ(d) \cong \Rep_d^{\nil} IQ (d) \times \mathbb{A}^1.$$ 
To prove (b), note that 
$$\Phi([\Rep^0_d\IQ])=\iota^\snil_*([\Rep^0_d\IQ\times\A^1])=\iota^\scal_*([\Rep^\scal_d\IQ])$$
and that $\Phi$ is $S_Q$-equivariant. By Proposition~\ref{prop:red-class} we have
	\begin{align*}
		\Phi\big(c\cap[\Rep^0_d\IQ]\big)
		&=\iota^\snil_*\big(c\cap[\Rep^\scal_d\IQ]\big)\\
		&=c\cap(\iota^\scal)_*[\Rep^\scal_d\IQ]\\
		&=c\cap[\Rep_d\Pi_Q]^\rd
	\end{align*}
The proof of (c) is more delicate. Consider the action on $\tRep_d \Pi_Q$ 
of the torus $\T=\GM$ which acts with weight $1$ on $\theta$, weight 
$-1$ on $x^*$ and weight $0$ on $x$. This action preserves the potential $W$ in \eqref{potential}.
Thus it commutes with the 
dimensional reduction isomorphisms \eqref{eq:star}. All of the above spaces, algebras and maps admit $\T$-equivariant versions 
$\bfA^{\nil,\T}_{IQ}$, etc., replacing everywhere $\H_*(-,\Q)$ by $\H^{\T}_*(-,\Q)$. These are modules over $\H^*_\T=\Q[t]$. 
The isomorphism 
of Lemma~\ref{lem:translation} is $\T$-equivariant, with $\T$ acting on $c\in\A^1$ with weight $1$. 
The self-intersection formula for the zero 
section $s: \Rep^\nil_dIQ \to \Rep^\snil_d IQ$ of the $\T$-equivariant line bundle 
$$\mathrm{pr}:\Rep^\snil_d IQ\cong\Rep^\nil_dIQ\times\A^1\to\Rep_d^\nil\IQ,$$
yields 
$$(s^\T)_*(x)=t\cdot(\mathrm{pr}^\T)^*(x).$$
Hence 
$$\iota^{\nil,\T}_*=t\cdot\Phi^\T:\bfA^{\nil,\T}_{IQ}(d)\to\bfA^\T_{IQ}(d),\quad (d\neq 0).$$
		
For $x\in\bfA^{\nil,\T}_{IQ}(d)$ and $x'\in\bfA^{\nil,\T}_{IQ}(d')$ with $d,d'\neq 0$, we have
$$t\,\Phi^\T(x*x')=\iota^{\nil,\T}_*(x*x')=\iota^{\nil,\T}_*(x)\,\iota^{\nil,\T}_*(x')=t^2\,\Phi^\T(x)*\Phi^\T(x').$$
From the purity in Proposition~\ref{prop: purity} we deduce that
$\bfA_{IQ}^\T \cong \bfA_{\Pi_Q}^\T$ is a free $\Q[t]$-module with $\bfA_{IQ}^\T/t\bfA_{IQ}^\T=\bfA_{IQ}$. Hence 
$$\Phi^\T(x*x')=t\,\Phi^\T(x)*\Phi^\T(x').$$
Reducing modulo $t$, i.e., applying the forgetful maps,
which are compatible with pushforward, smooth pullback and products, gives 
$$\Phi(\bar x*\bar x')=0.$$ 
To conclude, we use the fact that each class of $\bfA^\nil_{IQ}(d)$ lifts to $\bfA^{\nil,\T}_{IQ}(d)$ since $\bfA^\nil_{IQ}$ is pure.
\end{proof}

\medskip

\medskip\subsection{The generation theorem for $\bfA_{\Pi_Q}$}
The following proposition  is proved in the next section. 

\begin{proposition}\label{prop:nilp generation}
The nilpotent COHA $\bfA^\nil_{IQ}$ is generated by $\bigoplus_d \widehat{\g}^0(d)$ as an algebra.	
\end{proposition}

Using this result, we can prove the generation theorem Theorem~$\ref{thm:Markman2}$ for $\bfA_{\Pi_Q}$.

\begin{proof}[Proof of Theorem~$\ref{thm:Markman2}$] The theorem follows by combining 
Propositions~\ref{prop:seminilpcharacterizationBPS}, ~\ref{prop:Phi} and ~\ref{prop:nilp generation}. Indeed, we have
	$$\bfA^\nil_{IQ}(d)=\widehat{\g}^0(d)+\sum_{\substack{d'+d''=d,\\ d',d''\neq0}}\bfA^\nil_{IQ}(d')*\bfA^\nil_{IQ}(d'').$$
	Applying $\Phi$ we get
	$$
	\widehat{\g}^{IQ}(d)=\Phi(\bfA^\nil_{IQ}(d))=\Phi(\widehat{g}^0(d))=\widehat{\g}^{IQ}_{\taut}(d)
	$$
The theorem follows by \eqref{eq:star}.
\end{proof}

\medskip

\medskip\subsection{The generation theorem for $\bfA_{IQ}^\nil$} 
In order to prove Proposition~$\ref{prop:nilp generation}$, we first introduce a suitable stratification of $\Rep^\nil_dIQ$ by smooth and cohomologically pure pieces.

\smallskip\subsubsection{Jordan types}
A {Jordan type} of dimension $d$ is a tuple of dimension vectors 
$$\ud=(d^1,\ldots,d^s)
,\quad
d^i\in\N Q_0
,\quad
d^s\neq0
,\quad
|\ud|=\sum_i i\,d^i=d.$$
For such $\ud$ put
\[
e_i=d^i+d^{i+1}+\cdots+d^s
,\quad 
1\leqslant i\leqslant s
,\quad 
e_{s+1}=0,
\]
We have
$$e_1\geqslant e_2\geqslant\cdots\geqslant e_s>0
,\quad
\sum_ie_i=d
,\quad
d^i=e_i-e_{i+1}$$
and
\begin{equation}\label{eq:dimker thetai}
	e_1+\cdots+e_i=\dim\ker\theta^i.
\end{equation}
Let
\[
J(\ud)=\Big\{(x,\theta)\in\Rep^\nil_d\IQ\;;\;\theta^s=0\,,\, \dim\big(\ker\theta^i/\ker\theta^{i-1}\big)=e_i\ \text{for }1\leqslant i\leqslant s\Big\}.
\]
Thus $d^i$ is the dimension vector of the multiplicity space of Jordan blocks of size $i$. 
Since every $\ker\theta^i$ is a subrepresentation
$$\left(\Rep^\nil_d\IQ\right)^\cl=\bigsqcup_{|\ud|=d}J(\ud)$$ 
is a finite decomposition into locally closed substacks.
For $s=1$ we get 
$$J((d))=\Rep^0_d\IQ \cong \Rep_d Q.$$

\smallskip\subsubsection{Tautological bundles on $J(\ud)$}
On $J(\ud)$ the subspaces $\ker\theta^i$ form vector bundles, which are subrepresentations of the tautological quiver bundle. Multiplication by $\theta$ induces injections 
$$\ker\theta^{i+1}/\ker\theta^i\hookrightarrow\ker\theta^i/\ker\theta^{i-1}.$$ Let
\begin{align}\label{V}
	\calV^i=\coker\Big(\ker\theta^{i+1}/\ker\theta^{i}\xrightarrow{\ \theta\ }\ker\theta^{i}/\ker\theta^{i-1}\Big),
\end{align}
It is a $\C Q$-representation of dimension $e_i-e_{i+1}=d^i$. 

\begin{lemma}\label{lem:jordan}
The morphism
\begin{align}\label{pi}
	\pi_{\ud}\colon J(\ud)\to\prod_{i=1}^s\Rep_{d^i}Q,\quad(x,\theta)\mapsto(\calV^i)_i,
\end{align}
is an iterated vector bundle stack. Hence $J(\ud)$ is smooth and irreducible, the Borel-Moore homology space
$\H_*(J(\ud),\Q)$ is pure, and there is an isomorphism
$$\bigotimes_{i=1}^s\H^*_{GL_{d^i}} \cong \H_*(J(\ud),\Q), \quad
 \bigotimes_ic_i \mapsto \pi_{\ud}^*\Big(\bigotimes_i c_i\Big) \cap [J(\ud)] .$$
Moreover, we have
\begin{align}\label{dimJ}
	\dim J(\ud)=-\sum_{i,j}\min(i,j)\langle d^i,d^j\rangle=-\sum_{l=1}^s\langle e_l,e_l\rangle.
\end{align}
\end{lemma}

\begin{proof} This is a general result, valid for the tangent bundle to the stack of objects in a category of global dimension one
and not only for  $\Rep_d IQ$, see, e.g., \cite{MS20}.
%%%%%%%%%%%%%%%%%%%%%%%%%%%%%%
\end{proof}

\smallskip\subsubsection{The order and closures}
For $a,a' \in \N Q_0$ we write $a'<a$ if $a-a' \in \N Q_0 \smallsetminus \{0\}$. 
For Jordan types $\ud$ and $\ud'$ of the same dimension $d$, 
with associated sequences $(e_i)_i$ and $(e'_i)_i$, 
we write $\ud<\ud'$ if $\ud\neq\ud'$ and there exists $l$ such that for each $i<l$ we have
$$
e'_1+\cdots+e'_i=e_1+\cdots+e_i
,\quad 
e'_1+\cdots+e'_l\ <\ e_1+\cdots+e_l.
$$
By \eqref{eq:dimker thetai}, 
we have $\ud<\ud'$ if and only if $\dim\ker\theta^i$ is the same for $\ud$ and for $\ud'$ if $i<l$ and strictly larger for $\ud$ if $i=l$. 
Write $\ud\leqslant\ud'$ if either $\ud<\ud'$ or $\ud=\ud'$.
We consider the substacks 
$$Z_{\leqslant\ud}=\bigcup_{\ud'\leqslant\ud}J(\ud')
,\quad
Z_{<\ud}=\bigcup_{\ud'<\ud}J(\ud')$$

\begin{lemma}\label{lem:closed}
For every $\ud$ the substacks $Z_{\leqslant\ud}$ and $Z_{<\ud}$ are closed in $\Rep^\nil_d\IQ$.
\end{lemma}

\begin{proof}
Since there are finitely many strata, it is enough to consider the case of the substacks $Z_{\leqslant\ud}$. The functions $$(x,\theta)\mapsto\dim\ker\theta^i\in\N Q_0$$ are upper semicontinuous. If $y$ lies in the closure of the stratum of $x$, then 
$$\dim\ker\theta_y^i\geqslant\dim\ker\theta_x^i.$$ 
Hence at the first $i$ where they differ, the inequality is strict, i.e., the type of $y$ is $\leqslant$ the type of $x$. So $Z_{\leqslant\ud}$ is closed under specialization. It is constructible, being a finite union of locally closed strata, hence closed.
\end{proof}

\begin{corollary}\label{cor:purity-nil}
For every $\ud$ there is a short exact sequence
$$
0\to\H_*(Z_{<\ud},\Q)\to\H_*(Z_{\leqslant\ud},\Q)\to\H_*(J(\ud),\Q)\to0 .
$$
The mixed Hodge structure of $\bfA_{IQ}^\nil(d)$ is pure. 
\end{corollary}

\begin{proof}
Fix a total order $\preceq$ on Jordan types refining $\leqslant$. By Lemma~\ref{lem:closed}, the substacks 
$Z_{\preceq \ud}=\bigcup_{\ue \preceq \ud} J(\ue)$ are closed. Each step of the resulting filtration gives a localization long exact sequence, 
all of whose terms are pure by Lemma~\ref{lem:jordan} and induction. Hence, the long exact sequences split.
\end{proof}

\smallskip\subsubsection{Ordered products}~\label{sec:ordered}
Fix $\ud$ with associated sequence $\ue=(e_1,\ldots,e_s)$.
Set $\uf=(e_s,\ldots,e_1)$.
We consider the nilpotent COHA induction diagram for flags of subobjects with subquotients of dimensions $e_s,\ldots,e_1$:
\begin{equation}\label{diag:inductionnilp}
\prod_i\Rep^\nil_{e_i}\IQ\xleftarrow{\ q_\uf\ }\mathfrak{Filt}^\nil_\uf\IQ\xrightarrow{\ p_\uf\ }\Rep^\nil_d\IQ.
\end{equation}
We have
$$a_s*\cdots*a_1=(p_\uf)_*(q_\uf)^!(a_s\boxtimes\cdots\boxtimes a_1).$$
We define the stack
$$\mathfrak{Filt}^0_\uf\IQ=q_\uf^{-1}\Big(\prod_i\Rep^0_{e_i}\IQ\Big).$$ 
This stack parametrizes flags $0=F_0\subset F_1\subset\cdots\subset F_s$ of 
subrepresentations of a $d$-dimensional representation $(x,\theta)$ with $\theta(F_i)\subseteq F_{i-1}$ and $\dim F_i/F_{i-1}=e_i$. Observe 
that if $\theta(F_i)\subseteq F_{i-1}$ for all $i$, then $F_i\subseteq\ker\theta^i$. Hence $\dim\ker\theta^i\geqslant e_1+\cdots+e_i$ for all $i$, and 
the Jordan type of $(x,\theta)$ is $\leqslant\ud$. Thus, we have 
$$p_{\uf}(\mathfrak{Filt}^0_\uf\IQ) \subseteq Z_{\leqslant \ud}.$$
Let $q^0_\uf$ and $p^0_\uf$ be the restrictions of the maps $p_\uf,$ $ q_{\uf}$ to $\mathfrak{Filt}^0_\uf IQ$.
We  define
$$\Big(\mathfrak{Filt}^{0}_\uf\IQ\Big)^\circ=(p^0_\uf)^{-1}(J(\ud)),$$ 
with the open embedding 
$$j\colon\left(\mathfrak{Filt}^{0}_\uf\IQ\right)^\circ\hookrightarrow\mathfrak{Filt}^{0}_\uf\IQ .$$
The induction diagram \eqref{diag:inductionnilp} fits into the following commutative diagram
\begin{equation}\label{diag:inductionnilp2}
\begin{tikzcd}
\prod_i \Rep^\nil_{e_i}IQ  & \mathfrak{Filt}^\nil_\uf IQ \arrow[l,"q_{\uf}"'] & (\mathfrak{Filt}^\nil_\uf IQ)^\cl \arrow[l,"i"'] \arrow[r,"p_{\uf}"] & 
\Rep^\nil_{d}IQ\\
\prod_i \Rep^0_{e_i}IQ  \arrow[u,hook,"\iota"] & \mathfrak{Filt}^{0}_\uf IQ\arrow[u,hook,"\iota"] \arrow[l,"q^0_{\uf}"'] & (\mathfrak{Filt}^{0}_\uf IQ)^\cl \arrow[u,hook,"\iota"]\arrow[l,"i"'] \arrow[r,"p^0_{\uf}"] & Z_{\leqslant \ud}\arrow[u,hook,"\iota"]\\
	& (\mathfrak{Filt}^{0}_\uf IQ)^\circ \arrow[u,"j"] & (\mathfrak{Filt}^{0}_\uf IQ)^{\circ,\cl} \arrow[u, "j^\cl"] \arrow[l,"i^\circ"']\arrow[r,"p^{00}_{\uf}"]  & J(\ud) \arrow[u] .
\end{tikzcd}
\end{equation}
The pushforward functors $i_*$ and $ i^\circ_*$ are isomorphisms in Borel-Moore homology.
Note that the closed substack $Z_{\leqslant \ud}$ has no obvious derived enhancement. 
Finally, recall the bundles $\calV^i$ on $J(\ud)$ defined in \eqref{V}.
Let $\calE^i$ be the tautological bundle on $\Rep^0_{e_i}\IQ$. We set $\calE^{s+1}=0$. 

\begin{lemma}\label{lem:corresp}
The following hold:
\begin{enumerate}[label=$\mathrm{(\alph*)}$,leftmargin=8mm,itemsep=1mm]
	\item  The maps $i^\circ$ and $p^{00}_\uf\colon(\mathfrak{Filt}^{0}_\uf IQ)^{\circ,\cl}\to J(\ud)$ are isomorphisms.
	\item In $K$-theory we have
		$$
		(p^{00}_\uf)^*[\calV^i]=j^*(q^0_\uf)^*\big([\calE^i]-[\calE^{i+1}]\big)
		,\quad
		1\leqslant i\leqslant s.
		$$
	\item In Borel-Moore homology we have
	$$j^*(q^0_\uf)^!\big(\boxtimes_i[\Rep^0_{e_i}\IQ]\big)=[(\mathfrak{Filt}^{0}_\uf IQ)^\circ].$$
	\end{enumerate}
\end{lemma}

\begin{proof} The second statement of (a) follows from the fact that if $(x,\theta) \in J(\ud)$ then $F_i=\ker\theta^i$ is the unique flag of 
subrepresentations such that $\theta(F_i)\subset F_{i-1}$ and $\dim\;F_i/F_{i-1}=e_i$. Next, we argue that $(\mathfrak{Filt}^0_\uf IQ)^\circ$ is 
classical. As $j$ is an open embedding and $q^0_\uf$ is quasi-smooth, $(\mathfrak{Filt}^0_\uf IQ)^\circ$ is itself quasi-smooth. Moreover, the 
relative (virtual) dimension of $q_\uf$ is 
	$$-\sum_{i<j}(e_j,e_i)_{\C Q[\theta]}=0,$$ 
because the Euler form on $\C Q[\theta]$ vanishes. Thus we have
$$\vdim \; \mathfrak{Filt}^{0}_\uf IQ=\vdim\; \prod_i \Rep^0_{e_i}IQ=-\sum_i \langle e_i,e_i\rangle=\vdim\; J(\ud) $$
by \eqref{dimJ}. Thus the dimension of $(\mathfrak{Filt}^{0}_\uf IQ)^\circ$ coincides with that of its classical truncation. An irreducible quasi-
smooth derived stack whose classical truncation is smooth of dimension equal to the virtual 
dimension has a vanishing obstruction sheaf, hence is classical, and its virtual fundamental class is the fundamental class
by Theorem \ref{thm:mainB}. We turn to (b). If $(x,\theta,F_\bullet) \in (\mathfrak{Filt}^{0}_\uf IQ)^\circ$ then 
$F_i/F_{i-1}=\ker\theta^i/\ker\theta^{i-1}$, 
which is the pullback of $\calE^i$ under the map $q^0_{\uf}$. The claim follows from the short exact sequence 
$$0 \to \calE^{i+1} \stackrel{\theta}{\to} \calE^i \to \calV^i \to 0.$$
Finally, statement (c) is a property of quasi-smooth pullback.
\end{proof}

For all $i$, the algebra $\H^*(\Rep_{d^i}Q,\Q) \cong \H_{GL_{d^i}}^*$ is a free polynomial algebra in the Chern 
classes $\ch_k(\calV^i_j)$ of the tautological bundles $\calV^i_j$ for $j \in Q^0$ and $k=1, \ldots, d^i_j$. Thus
there is a well-defined algebra morphism
$$
h_\ud\colon\bigotimes_{i=1}^s\H_{GL_{d^i}}^*\to
\bigotimes_{i=1}^s\H^*(\Rep^0_{e_i}IQ,\Q),\quad \ch_k(\calV^i)\mapsto \ch_k(\calE^i)-\ch_k(\calE^{i+1}) .
$$
%where $c(\mathcal{U})=1+ c_1(\mathcal{U}) + c_2(\mathcal{U}) + \cdots$ stands for the total Chern class.
Recall the map $\pi_\ud$ in \eqref{pi}. By Lemma~\ref{lem:corresp}(b) we have 
\begin{equation}\label{eq:proofgeneration1}
(p^{00}_\uf)^*\pi_\ud^*(c)=j^*(q^0_\uf)^*h_\ud(c).
\end{equation}
For every $\ud$ and every $c\in\bigotimes_i\H_{GL_{d^i}}^*$ we define
$$
\mu_\ud(c)=(p_\uf)_*(q_\uf)^!\Big(h_\ud(c)\cap\boxtimes_i[\Rep^0_{e_i}\IQ]\Big)\ \in\ \sum\ \widehat{\g}^0(e_s)*\cdots*\widehat{\g}^0(e_1)
$$

\begin{lemma}\label{prop:jordan-product}\hfill
\begin{enumerate}[label=$\mathrm{(\alph*)}$,leftmargin=8mm,itemsep=1mm]
\item
For every $\ud$ and every $c\in\bigotimes_i\H_{GL_{d^i}}^*$ the class $\mu_\ud(c)$ is supported on $Z_{\leqslant\ud}$. 
\item
Its restriction to the open stratum $J(\ud)\subset Z_{\leqslant\ud}$ is $\pi^*_\ud(c)\cap[J(\ud)]$.
\end{enumerate}
\end{lemma}

\begin{proof}
By base change and the projection formula, we have
$$(q_{\uf})^!\left(h_\ud(c)\cap \boxtimes_i[\Rep^0_{e_i}\IQ]\right)=
\iota_*\Big((q_{\uf}^0)^!\left(h_\ud(c)\cap\boxtimes_i[\Rep^0_{e_i}\IQ]\right)\Big).$$ 
Hence the left hand side is supported on $\mathfrak{Filt}^0_{\uf}IQ$. 
Its pushforward by $p_\uf$ is therefore supported on 
$$\Im(p^0_\uf)\subseteq Z_{\leqslant\ud}.$$ 
This proves part (a). For (b), note that restriction to the open substack $J(\ud)$ commutes 
with proper pushforward via the right Cartesian square in \eqref{diag:inductionnilp2}. Using Lemma~\ref{lem:corresp}(c) 
and \eqref{eq:proofgeneration1}, we get
\begin{align*}
	j^*(p_{\uf})_*(q_\uf)^!\Big(h_\ud(c)\cap\boxtimes_i[\Rep^0_{e_i}\IQ]\Big)
	&=(p^{00}_\uf)_*\Big(j^*(q^0_\uf)^*h_\ud(c)\cap[\mathfrak{Filt}^0_\uf\IQ]\Big)\\
	&=(p^{00}_{\uf})_*\Big((p^{00}_\uf)^*(\pi_\ud)^* (c)\cap[\mathfrak{Filt}^0_\uf\IQ]\Big)\\
	&=\pi_\ud^*(c)\cap[J(\ud)] .
\end{align*}
Finally, the class $h_\ud(c)$ being a sum of pure tensors 
$\otimes_ia_i$ with
$$a_i\in\H^*(\Rep^0_{e_i}IQ,\Q)=\H^*_{GL_{e_i}},$$ 
the class $\mu_\ud(c)$ is a sum of ordered products of elements of $\widehat{\mathfrak g}^0$ of the form
$$(a_s\cap[\Rep^0_{e_s}\IQ])*\cdots*(a_1\cap[\Rep^0_{e_1}\IQ]) .$$
\end{proof}

\smallskip\subsubsection{Proof of Proposition~$\ref{prop:nilp generation}$.} We are now ready to finish the proof of Markman's theorem. Fix a 
total order $\preceq$ on the set of Jordan types of dimension $d$ refining $\leqslant$.
Put $Z_{\preceq \ud}=\bigcup_{\ue \preceq \ud}J(\ue)$. 
These are closed substacks by the proof of Lemma~\ref{lem:closed}.
%and we have $Z_{\preceq\ud_k}\subseteq Z_k$. 
By Corollary~\ref{cor:purity-nil}, the pushforwards by the inclusion 
$$Z_{\preceq \ud} \to \Rep^\nil_d IQ$$ is injective.
This yields a filtration 
$$\H_*(Z_{\preceq \ud}\,,\,\Q)\subseteq \bfA_{IQ}^\nil(d)$$ 
whose associated graded is given by
$$\H_*(Z_{\preceq \ud}\,,\,\Q)/\H_*(Z_{\prec \ud}\,,\,\Q) \cong\H_*(J(\ud)\,,\,\Q).$$ 
By Lemma~\ref{prop:jordan-product}, we have
$$\mu_{\ud}(c) \in (\pi_{\ud}(c)\cap [J(\ud)]) \;\oplus\; \H_*(Z_{\prec \ud}\,,\,\Q), \quad \forall c \in \bigotimes_{i} \H^*_{GL_{d_i}} .$$
By induction, this implies that the classes $\mu_{\ud}(c)$ for all $\ud$ 
and $c$ running over a basis of $\bigotimes_{i} \H^*_{GL_{d_i}}$  form a basis of $\H_*(\Rep^\nil_dIQ\,,\,\Q)$. 
We deduce that for every $d$ we have
$$
\bfA_{IQ}^\nil(d)=\sum_{|\ud|=d}\ \sum \widehat{\g}^0(e_s)*\cdots*\widehat{\g}^0(e_1),
$$
where the inner sum runs over all Jordan types of weight $d$ and $\ue$ is the sequence associated with $\ud$. 
This proves Proposition~\ref{prop:nilp generation} and Theorem~\ref{thm:Markman2}.

%%%%%%%%%%%%%%%%%%%%%%%%%%%%%%%%%%%%%%%%%%%%%%%%%%%%%

%=====================================================================
\medskip\section{Character varieties and Higgs moduli.}
%=====================================================================

\medskip

In this section we fix the notation and briefly collect known results on the (co)homology of twisted character varieties and the (co)homology of 
moduli stacks of Higgs bundles on smooth projective curves. In particular, we recall the  $C=W$ and P=C theorems as well as 
the curious Hard Lefschetz property. We refer to \cite{HRV08}, \cite{S17}, \cite{M19}, \cite{MS} and \cite{HMMS} for details.
Throughout, we fix an integer $g >1$, a compact topological surface $\Sigma_g$ of genus $g$, 
and a smooth complex projective curve $C$ of genus $g$.

\medskip\subsection{Character varieties}  Let us fix a pair of coprime integers $(r,d)$ with $r>0$ and define the genus $g$ twisted character variety $M^B(r,d)$ by \eqref{eq:defcharacter}. It is a smooth affine algebraic variety of dimension $2N$. Its cohomology 
$\H^*(M^B(r,d),\Q)$ carries a mixed Hodge structure which is of Hodge-Tate type. The weight filtration is
$$W_{\ell}\H^m(M^B(r,d),\Q), \quad \ell \in [m,2m] .$$

We next describe tautological generators for the cohomology ring $\H^*(M^B(r,d),\Q)$. 
The quotient $\tilde{M}^B(r,d)$ of $M^B(r,d)$ by the action of $(\GM)^{2g}$ given by rescaling the 
generators $(A_i, B_i)_{i=1, \ldots g}$ is a connected component of the $PGL_r(\C)$ (untwisted) character variety.
There is a ring isomorphism
\begin{equation}\label{eq:cohcharPGL}
	\H^*(M^B(r,d), \Q) \cong \H^*(\tilde{M}^B(r,d)) \otimes \H^*((\C^*)^{2g} ,\Q).
\end{equation}
There is a tautological (non-algebraic) $PGL_r(\C)$-principal bundle $\tilde{\mathcal{V}}$ on $\tilde{M}^B(r,d) \times \Sigma_g$ yielding 
tautological classes
$$\overline{c}_k(\eta)=\int_{\Sigma_g} \overline{c}_k(\tilde{\mathcal{V}}) \cup p_2^*(\eta) \in \H^{2k-2+|\eta|}(\tilde{M}^B(r,d),\Q)
, \quad 
\forall k \geqslant 2, \eta \in \H^*(\Sigma_g,\Q)$$
where $p_2: \tilde{M}^B(r,d) \times \Sigma_g \to \Sigma_g$ is the obvious projection and $\overline{c}_k$ is given, in terms of Chern roots $\{\zeta_i\}_i$ by the $k$th elementary symmetric function of $\{\zeta_i-\overline{\zeta}\}_i$ where $\overline{\zeta}$ is the average of the $\zeta_i$. These classes are denoted
$$\{\alpha_k, \beta_k, \phi_{i,k}\},\quad
i =1, \ldots, 2g, \quad
k=2, \ldots, r$$ and are defined as
$$\alpha_k=\overline{c}_k(1) \in \H^{2k-2}(\tilde{M}^B(r,d),\Q), \quad \beta_k=\overline{c}_k([\Sigma_{g}]) \in \H^{2k}(\tilde{M}^B(r,d),\Q),$$
$$ \phi_{i,k}=\overline{c}_k(\gamma_i) \in \H^{2k-1}(\tilde{M}^B(r,d),\Q)$$
where $\{\gamma_1, \ldots, \gamma_{2g}\}$ form a basis of $\H^1(\Sigma_g,\Q)$. 
By Markman's theorem \cite[Cor.~2]{Markman}, the classes $\{\alpha_k, \beta_k, \phi_{i,k}\}$ generate the cohomology ring 
$\H^*(\tilde{M}^B(r,d),\Q)$, yielding in turn a collection of generators for $\H^*(M^B(r,d),\Q)$ by \eqref{eq:cohcharPGL}. 
As shown by Hausel and Rodriguez-Villegas, the pure part of $\H^*(M^B(r,d),\Q)$ is the subring generated by 
$\{\beta_k\}$, see \cite[Cor.4.1.10]{HRV08}, 
\begin{equation}\label{eq:purepartBetti}
\H^{*}_{\text{pure}}(M^B(r,d),\Q)=\bigoplus_m W_m\H^m(M^B(r,d),\Q)=\Q\langle\beta_2, \beta_3, \ldots,\beta_r\rangle .
\end{equation} 
More generally, the place in the weight filtration of each of the above generators
is determined by Shende \cite{S17} who shows that the Chern degree coincides with the weight (C=W):
\begin{equation}\label{eq:shendeweights}
\overline{c}_k(\eta) \in W_{2k}\H^{2k-2+|\eta|}(\tilde{M}^B(r,d),\Q) , \quad k \geqslant 2, \eta \in \H^*(\Sigma_g,\Q) ,
\end{equation}
and the degree $1$ classes generating $\H^*(\GM,\Q)^{2g}$ all have weight $2$.

\medskip\subsection{Moduli spaces of Higgs bundles} 
Fix a pair $(r,d)$ of coprime integers with $r>0$. Let $M^D(r,d)$ be the moduli space of semistable Higgs bundles
$$(\mathcal{E},\varphi), \quad\varphi \in \text{Hom}_{\mathcal{O}_C}(\mathcal{E},\mathcal{E}\otimes \Omega_C)$$
of rank $r$ and degree $d$ on $C$. Thanks to the coprimality assumption on $(r,d)$, semistable Higgs bundles of rank $r$ and degree $d$ 
are automatically stable, and $M^D(r,d)$ is a smooth quasi-projective algebraic variety, equipped with a natural symplectic structure. It is of 
dimension $2N$ with
$$N=(g-1)r^2+1.$$ 
Scaling the Higgs field $\varphi$ defines a $\GM$-action, which endows $M^D(r,d)$ with the structure of a 
semi-projective variety, see \cite{HRVsemiprojective}. In particular, its cohomology carries a pure mixed Hodge structure. 
Let $$h: M^D(r,d) \to \mathcal{A}_r=\bigoplus_{k=1}^r \H^0(C,\Omega_C^{\otimes k})\cong \mathbb{A}^N.$$
denote the Hitchin map. It is a Lagrangian fibration.
In particular, it is proper. 
The Hitchin map induces a perverse filtration on $\H^*(M^D(r,d),\Q)$ such that
\begin{align}\label{PF1}
\P_{k}\,\H^*(M^D(r,d),\Q)=\H^*\left(\calA_r,{}^p\tau_{\leqslant k+N} h_*(\Q_{M^D(r,d)})\right).
\end{align}
It is normalized to be nontrivial in the interval $[0,2N]$.

As for the character varieties, the cohomology of $M^D(r,d)$ is tautologically generated. Consider the action of $T^*Pic^0(C)$ on 
$M^D(r,d)$ given by 
$$(\mathcal{L}, \theta) \cdot (\calE, \varphi)=(\mathcal{L} \otimes \calE, \text{id}\otimes\varphi + \theta \otimes \text{id} ).$$ 
The quotient of $M^D(r,d)$ by this action is a connected component  $\tilde{M}^D(r,d)$ of the moduli space of stable $PGL_r$-Higgs bundles on $C$, and there is a canonical algebra isomorphism
\begin{equation}\label{eq:cohoHiggsPGL}
\H^*(M^D(r,d),\Q) \cong \H^*(\tilde{M}^D(r,d),\Q) \otimes \H^*(T^*Pic^0(C),\Q) .
\end{equation}
The scheme $\tilde{M}^D(r,d) \times C$ carries a tautological algebraic $PGL_r(\C)$-principal bundle 
$\tilde{\mathcal{V}}$. 
We define the tautological classes to be
$$\overline{c}_k(\eta)=\int_{C} \overline{c}_k(\tilde{\mathcal{V}}) \cup p_2^*(\eta) \in \H^{2k-2+|\eta|}(\tilde{M}^D(r,d),\Q)
, \quad 
\forall k \geqslant 2, \eta \in \H^*(C,\Q)$$
where $p_2: \tilde{M}^D(r,d) \times C \to C$ is the projection. 
By Markman's theorem again, the cohomology algebra $\H^*(\tilde{M}^D(r,d),\Q)$ is generated by the tautological classes 
$\overline{c}_k(\eta)$ for $k=2, \ldots, r$ and $\eta \in \H^*(C,\Q)$. 
Moreover, thanks to \cite{MS}, \cite{HMMS}, their perverse degree is equal to their Chern degree, i.e., we have
\begin{equation}\label{eq:P=CPGL}
\overline{c}_k(\eta) \in \P_k\,\H^{2k-2+|\eta|}(\tilde{M}^D(r,d),\Q), \quad \forall k \geqslant 2, \eta \in \H^*(C,\Q).
\end{equation}
The perverse degrees of the degree $1$ classes generating $\H^*(T^*Pic^0(C),\Q)$ are all equal to $1$.
The perverse filtration on $\H^*(M^D(r,d),\Q)$ is multiplicative (P=C). 
In particular, we have
\begin{equation}\label{eq:PurePerversecoprime}
	\bigoplus_k \P_k\,\H^{2k}(M^D(r,d),\Q) =\Q\langle \overline{c}_2([\text{pt}]), \overline{c}_3([\text{pt}]), \ldots, \overline{c}_r([\text{pt}]) \rangle . 
\end{equation}

\medskip

\medskip\subsection{Non-abelian Hodge correspondence and P=W} 
Fix an isomorphism of topological spaces $\Sigma_g \cong C$.
This identifies $\H^*(\Sigma_g,\Q)$ with $\H^*(C,\Q)$. 
The non-abelian Hodge correspondence of Simpson and Corlette provides, for any coprime pair $(r,d)$, a diffeomorphism
$$\Theta: M^B(r,d) \stackrel{\cong}{\longrightarrow} M^D(r,d)$$
which is compatible with tautological classes in the following sense: the associated ring isomorphism 
$$\Theta_*: \H^*(M^B(r,d),\Q) \stackrel{\cong}{\longrightarrow} \H^*(M^D(r,d),\Q)$$
maps $\overline{c}_k(\eta)$ to $\overline{c}_k(\eta)$ for all $k \geqslant 2, \eta \in \H^*(C,\Q)$ and restricts to the isomorphism $$\H^*(\GM,\Q)^{2g} \cong \H^*(T^*Pic^0(C),\Q)$$
induced by the homeomorphism $\Sigma_g \cong C$. The P=W identity is the coincidence between \eqref{eq:shendeweights} and \eqref{eq:P=CPGL}. 

\medskip\subsection{Curious Hard Lefschetz} There is an extra symmetry on $\H^*(M^D(r,d),\Q)$, provided on the Higgs moduli side by the relative Hard Lefschetz isomorphism. 
For $\alpha \in \H^2(M^D(r,d),\Q)$ a class which is relatively ample with respect to $h$, capping with $\alpha$ yields isomorphisms
$$
\alpha^{j}\colon \gr^P_{N-j}\,\H^{m}(M^D_{r,d},\Q)\xrightarrow{\cong}\gr^P_{N+j}\,\H^{m+2j}(M^D_{r,d},\Q)
$$
for any $j,m$. Taking $j=N-k$ and $m=2k$, and using the fact that 
$$\gr^P_{k}\,\H^{2k}(M^D_{r,d},\Q)=\P_{k}\H^{2k}(M^D_{r,d},\Q),$$ 
e.g., by P=W and the bound $W_m\H^\ell=0$ for $m \not\in [\ell, 2\ell]$, this gives 
\begin{equation}\label{eq:CHLpure}
\P_{k}\,\H^{2k}(M^D(r,d),\Q)\cong \gr^P_{2N-k}\,\H^{2N}(M^D(r,d),\Q).
\end{equation}
Since $M^D(r,d)$ retracts to the global nilpotent cone, we have
$$j>2N\Rightarrow \H^{j}(M^{D}(r,d),\Q)=0.$$
Hence, the relative Hard Lefschetz isomorphism implies that
$$P_{N-1}\H^{2N}(M^D(r,d),\Q)=0 . $$
In other words, \eqref{eq:CHLpure} relates the perverse pure part of $\H^*(M^D(r,d),\Q)$ with the perverse filtration of the 
{middle} cohomology of $M^{D}(r,d)$. 
Thus, using \cite[Cor.~1.7]{S}, \cite{M20} and \cite[Rem. 4.4.6]{HRV08}, we deduce that
\begin{equation}\label{eq:lowerbound}
\dim\;\Big(\bigoplus_k \P_k\,\H^{2k}(M^D(r,d),\Q)\Big)=\dim\; H^{2N}(M^D(r,d),\Q)=A_{g,r}(0)=A_{S_g,r}(1) ,
\end{equation}
where $A_{g,r}(\underline{\sigma})$ is the Kac polynomial counting indecomposable vector bundles of rank $r$. 
Note that this dimension is also the number of irreducible components of the semistable global nilpotent cone in $M^D(r,d)$,
see \cite[Cor.1.9]{S}.

\medskip

%=====================================================================
\medskip\section{The Higgs COHAs and the BPS Lie algebra $\widehat{\mathfrak{g}}^C$. }
%=====================================================================

In this section, we drop the coprimality assumption on $(r,d)$ and introduce the BPS Lie algebra of the stack of semistable Higgs sheaves 
on $C$. We recall the $\chi$-independence theorem of Kinjo-Koseki as well the  stacky P=C conjecture of Picovschi. We refer to 
\cite{KinjoKoseki}, \cite{DHKSV} and \cite{Raphael} for details.

\medskip

\medskip\subsection{The Higgs COHA} The slope of a coherent sheaf $\mathcal{F}$ of class $\alpha=(r,d)$, i.e., 
of rank $r$ and degree $d$, is
$$\mu(\mathcal{F})=\mu(\alpha)=\frac{d}{r} .$$
For $\alpha, \alpha' \in \Z^2$ we write
$$(\alpha,\alpha')=2(1-g)\,\text{rank}(\alpha)\,\text{rank}(\alpha').$$
Note that, this coincides with the Euler form \eqref{eq:Eulerform} for the $S_g$ quiver, 
which is the only quiver that we will consider in the rest of the paper.
Let $\mathfrak{M}^{\ss}(\alpha)$ be  the stack of semistable Higgs bundles on $C$ of class $\alpha$. 
It is a quasi-smooth, quasi-compact derived stack of virtual dimension
$$\vdim\;\mathfrak{M}^\ss(\alpha)=-(\alpha,\alpha) .$$
Fix a slope $\nu \in \Q$. The COHA of semistable Higgs bundles of slope $\nu$ is 
$$\bfA_{C,l}(\alpha)=\H_{l-(\alpha,\alpha)}(\mathfrak{M}^\ss(\alpha),\Q), \quad \bfA_C(\alpha)=\bigoplus_l \bfA_{C,l}(\alpha), \quad \bfA^{(\nu)}_C=\bigoplus_{\mu(\alpha)=\nu} \bfA_C(\alpha) .$$
There is an induction diagram in all ways similar to \eqref{diag:defprodcoha1} giving rise to a $\N \times \Z$-graded associative algebra structure on $\bfA_C^{(\nu)}$.

\medskip

\medskip\subsection{The BPS Lie algebra $\widehat{\g}^C$} When $\alpha=(r,d)$ is primitive, i.e., for coprime $(r,d)$, 
the classical truncation $\mathfrak{M}^{\ss,\cl}(\alpha)$ is a $\GM$-gerbe over the moduli space $M^D(r,d)$. 
In general there is a good moduli space
$$p: \mathfrak{M}^{\ss,\red}(r,d) \to M^D(r,d) .$$
For $\mu(\alpha)=\mu(\alpha')$ there is an addition map
$$\oplus_{\alpha,\alpha'}: M^D(\alpha) \times M^D(\alpha') \to M^D(\alpha+\alpha').$$
Let us fix a slope $\nu=\mu(\alpha_0)$ with $\alpha_0=(r_0,d_0)$ primitive. 
There is a collection of pure mixed Hodge modules $\mathcal{BPS}^C_{n\alpha_0}$ on $M^D(n\alpha_0)$ for $n \geqslant 1$ such that
\begin{equation}\label{eq:integrality2}
	\text{Sym}_{\oplus}\Big( \bigoplus_{n \geqslant 1} \mathcal{BPS}_{n\alpha_0}^C \otimes\H_{\GM}^*\Big) 
	\cong \bigoplus_{n \geqslant 1} p_*\Big(\DQ_{\mathfrak{M}^{\ss}(n\alpha_0)}\otimes \mathbb{L}^{-\frac{n^2}{2}(\alpha_0,\alpha_0)}\Big) .
\end{equation}
The COHA of semistable Higgs bundles enjoys the same structural properties as its quiver cousin, and, in fact,
as  the COHA of any $2CY$ category in the sense of \cite{DHSM2}. Namely, we define the BPS Lie algebras and affinized 
BPS Lie algebras as
\begin{align*}
	\g^{C,(\nu)}&=\bigoplus_{n \geqslant 1} \g^{C}(n\alpha_0)= \bigoplus_{n \geqslant 1} \H^{-*}\Big(M^D(n\alpha_0)\,,\,\mathcal{BPS}_{n\alpha_0}^C\Big),\\
	\widehat{\g}^{C,(\nu)}&=\bigoplus_{n \geqslant 1} \widehat{\g}^{C}(n\alpha_0)= \bigoplus_{n \geqslant 1} \H^{-*}\Big(M^D(n\alpha_0)\,,\,\mathcal{BPS}_{n\alpha_0}^C\otimes H_{\GM}^*\Big).
\end{align*}
They are graded subspaces in $\bfA^{(\nu)}_{C}$, where the generator $u \in \H_{\GM}^2$ is in degree $-2$, 
and both form Lie subalgebras under the COHA multiplication. 
Note that we only consider the semistable BPS Lie algebras, 
and not the {global} BPS Lie algebra considered in \cite{DHKSV}.
Morever, we have
\begin{equation}\label{eq:A=Ug2}
	\text{U}\big(\widehat{\g}^{C,(\nu)}\big) = \bfA_{C}^{(\nu)}
\end{equation}
and passing to global sections in \eqref{eq:integrality2} recovers the PBW isomorphism
\begin{equation}\label{eq:PBW2}
	\Sym \Big( \bigoplus_{n\geqslant 1} \widehat{\g}^{C}(n\alpha_0)\Big) \stackrel{\cong}{\longrightarrow} \bfA_{C}^{(\nu)}
\end{equation}
where the map is given by the COHA multiplication. Finally, there is a comultiplication \cite{DHKSV}
\begin{equation}\label{eq:coproduct2}
	\Delta: \bfA_{C}^{(\nu)} \to \bfA_{C}^{(\nu)} \otimes \bfA_{C}^{(\nu)}
\end{equation}
It is compatible with the COHA product.
The space of primitive elements is $\widehat{\g}^{C,(\nu)}$. 
The graded dimension of $\g^{C,(\nu)}$ is given in terms of
the Kac polynomial $A_{g,r}(\underline{\sigma})$ of \cite{S}.
We will not need this. 

\begin{remark}
The non-abelian Hodge correspondence yields an isomorphism of COHAs and BPS Lie algebras, involving $\bfA_C^{(0)}$ 
on the Higgs side and the COHA of untwisted local systems on $\Sigma_g$ on the Betti side, see \cite{HennecartNAHT}, \cite{PortaSala}. 
\end{remark}

\medskip

\medskip\subsection{The perverse filtration on the BPS Lie algebra}
The Hitchin map
$$h: M^D(r,d) \to \mathcal{A}_r$$
is well-defined for every $(r,d)$. It is proper and induces a perverse filtration on $\g^C(r,d)$ defined by
\begin{equation}\label{PF2}
\P_k\,\g_\ell^C(r,d)=\text{Im}\left(\H^{-\ell}\big(\mathcal{A}_r,{}^p\tau_{\leqslant k-N}h_*\mathcal{BPS}^C_{(r,d)}\big) \to \H^{-\ell}\big(\mathcal{A}_r,h_*\mathcal{BPS}^C_{(r,d)}\big)\right)
\end{equation}
This definition ensures that the perverse filtration on $\g^C(r,d)$ is non trivial only in the interval $[0,2N]$, i.e.,
we have $$P_{-1}\,\g^C(r,d) =\{0\}, \quad P_{2N}\,\g^C(r,d)=\g^C(r,d) .$$
For coprime $(r,d)$, every semistable Higgs bundle is stable and $M^D(r,d)$ is smooth of complex dimension $2N$.
In this case the BPS sheaf is the intersection complex
$$
\mathcal{BPS}_{(r,d)}^C=\IC(M^D(r,d),\Q)=\Q_{M^D(r,d)}[2N]
$$
Taking hypercohomology in degree $-i$ therefore gives
\begin{equation}\label{eq:PoincaréBPS}
\g_{i}^C(r,d) \cong \H^{2N-i}\left(M^D(r,d),\Q\right) .
\end{equation}
This is the form of Poincar\'e duality compatible with the Borel--Moore homology realization of the BPS Lie algebra:
the shift $[2N]$ changes BPS degree $-i$ into ordinary cohomological degree $2N-i$.
It also explains the normalization of the perverse filtration. Indeed,
$$
{}^p\tau_{\leqslant k-N}\big(h_*\Q_{M^D(r,d)}[2N]\big)
\cong \big({}^p\tau_{\leqslant k+N}h_*\Q_{M^D(r,d)}\big)[2N],
$$
so that the filtration \eqref{PF2}, under the identification \eqref{eq:PoincaréBPS}, is precisely \eqref{PF1}.

Definition \eqref{PF2} extends to the affinized BPS Lie algebra $\widehat{\g}^C$, by replacing
$\mathcal{BPS}^C_{(r,d)}$ with $\mathcal{BPS}^C_{(r,d)}\otimes\H^*_{\GM}$.
Write $\H^*_{\GM}=\Q[u]$, where $u$ has homological degree $-2$ and perverse degree $2$.
The canonical inclusion of the constant term induces an isomorphism
\begin{equation}\label{eq:isomgradedpureperverse}
\bigoplus_\ell \gr^\P_{N-\ell}\,\g^C_{2\ell}(r,d) \cong \bigoplus_\ell \gr^\P_{N-\ell}\,\widehat{\g}^C_{2\ell}(r,d) .
\end{equation}
Here the absence of contributions from positive powers of $u$ uses, in addition to their bidegrees,
the following vanishing:
$$
\P_{N-s-1}\,\g^C_{2s}(r,d)=0 \qquad (s\in\Z).
$$
For coprime $(r,d)$ this follows from \eqref{eq:PoincaréBPS}, P=W and the inequality
$W_{2k-1}\H^{2k}(M^B(r,d),\Q)=0$; for arbitrary $d$ it follows by the filtered
$\chi$-independence isomorphisms stated below in \eqref{eq:chi_indepBPS}, after choosing $d'$ coprime to $r$.
Indeed, keeping track of both gradings gives
$$
\gr^\P_{N-\ell}\,\widehat{\g}^C_{2\ell}(r,d)
\cong \bigoplus_{m\geqslant 0} u^m\,
\gr^\P_{N-\ell-2m}\,\g^C_{2(\ell+m)}(r,d).
$$
For $m\geqslant 1$ we have $N-\ell-2m<N-(\ell+m)$, so the corresponding summand vanishes
by the preceding bound with $s=\ell+m$. Only $m=0$ survives, which proves \eqref{eq:isomgradedpureperverse}.

\medskip\subsection{Tautological classes} Let $\mathcal{V}$ be the tautological vector bundle on $\mathfrak{M}^{\ss}(\alpha) \times C$. As before, we define tautological cohomology classes by
$$\ch_k(\eta)=\int_C \ch_k(\calV)\cup p_2^*(\eta) \in \H^{2k-2+|\eta|}(\mathfrak{M}^\ss(\alpha),\Q), \quad 
\forall k \geqslant 1, \eta \in \H^*(C,\Q)$$
and we define an algebra
$$S_C=\Q\left[\ch_k([\text{pt}]), \ch_k(1), \ch_k(\gamma_i)\;;\; k \geqslant 1, i=1,\ldots, 2g\right]$$
acting on $\bfA_C^{(\nu)}$ for any slope $\nu$. Lemma~\ref{lem:hopf-module} holds in the present setting as well, and thus the action of $S_C$ on $\bfA_C^{(\nu)}$ preserves $\widehat{\g}^{C,(\nu)}$.
We put 
$$u=c_1([\text{pt}]) \in \H^2(\mathfrak{M}^\ss(\alpha),\Q) .$$
Note that by definition, $u=c_1(\calV_{x})$ for any $x \in C$, where $\calV_x$ is the restriction of $\calV$ 
to $\mathfrak{M}^\ss(r,d) \times \{x\}$.

If $(r,d)$ is primitive, then the map
$$p:\mathfrak{M}^{\ss,\red}(r,d) \to M^D(r,d)$$ is a $\GM$-gerbe and we may consider a
$PGL_r$-principal bundle $\tilde{\mathcal{V}}$ on the moduli space $M^D(r,d)$. 
It is  related to the tautological bundle $\calV$ by the relation
$$(p\times \text{id})^*(\tilde{\mathcal{V}})=\calV \otimes (\det \calV)^{-1/r}$$
where $\det \calV$ is the universal line bundle pulled back from the map 
$$\det: \mathfrak{M}^{\ss,\red}(r,d) \to T^*Pic^d(C), \quad (\calE,\varphi) \mapsto (\det\calE, \tr(\varphi)) .$$ 
In particular, the tautological classes on $\mathfrak{M}^\ss(r,d)$ and $M^D(r,d)$ are related via the equality
$$\ch(\calV)=(p \times \text{id})^*\ch(\tilde{\mathcal{V}})\,e^{\frac{1}{r}c_1(\calV)}.$$
Note that $c_1(1)=d$,
the classes $\phi_{i,1}=c_1(\gamma_i)$ generate a copy of $\H^*(T^*Pic^d,\Q)$,
 and $u$ has a non-zero contribution in $\H^2_{\GM}$ for any trivialization of the $\GM$-gerbe $p$.
We deduce that
$$\H^*(\mathfrak{M}^\ss(r,d),\Q) = p^*(\H^*(M^D(r,d)))[u].$$

\medskip\subsection{The reduced virtual fundamental class and Markman's theorem for $\widehat{\g}^{C}$} 
We refer to \cite[Sec.7.2]{DHKSV} for the definition of the reduced derived stack $\mathfrak{M}^{\ss,\rd}(\alpha)$.
It is a closed substack of $\mathfrak{M}^\ss(\alpha)$. 
It should not to be confused with the reduced stack $\mathfrak{M}^{\ss,\red}(\alpha)$, which is a classical stack, 
and is a closed substack of $\mathfrak{M}^{\ss,\rd}(\alpha)$. 
The reduced derived stack is quasi-smooth and satisfies
$$(\mathfrak{M}^{\ss,\rd}(\alpha))^\cl=\mathfrak{M}^{\ss,\cl}(\alpha)
,\quad
\vdim(\mathfrak{M}^{\ss,\rd})(\alpha)=\vdim\mathfrak{M}^{\ss}(\alpha)+1=1-(\alpha,\alpha) .
$$
It is obtained from $\mathfrak{M}^{\ss}(\alpha)$ by removing the trivial one-dimensional summand in the obstruction space via the trace map 
$$\tr_{\mathcal{F}}:\text{Ext}^2(\mathcal{F},\mathcal{F})\cong  \text{Hom}(\mathcal{F},\mathcal{F})^\vee \to \C.$$
 %When $g>1$, $\mathfrak{M}^{\ss,\rd}(\alpha)$ coincides with the classical truncation $\mathfrak{M}^{\ss,\cl}(\alpha)$ whenever $\text{rank}(\alpha)>0$; indeed in that case $\mathfrak{M}^{\ss,\cl}(\alpha)$ is generically a $\GM$-gerbe over $M^D(\alpha)$ hence $\vdim\;\mathfrak{M}^{\ss,\cl}(\alpha)=\vdim\;\mathfrak{M}^{\ss,\rd}(\alpha)$ and we may apply \cite[Prop.1]{Khan}.
We define
$$\widehat{\g}^{C}_{\taut}(\alpha)=S_C \bullet [\mathfrak{M}^{\ss,\rd}(\alpha)] .$$
The following may be viewed as an analog of Markman's generation theorem for the cohomology of $M^D(r,d)$ 
for nonprimitive pairs $(r,d)$. 

\begin{theorem}[{\cite[Thm.7.5]{DHKSV}}]
For all classes $\alpha$ of positive rank we have
$$\widehat{\g}^{C}_{\taut}(\alpha)=\widehat{\g}^{C}(\alpha) .$$
\end{theorem}

\medskip\subsection{$\chi$-independence} 
 The Hitchin base $\mathcal{A}_r$ is independent of $d$. 
 Fixing $r$, it thus makes sense to compare the semisimple 
 constructible complexes $h_*(\mathcal{BPS}^C({r,d})) \in D^b(\mathcal{A}_r)$ for various values of $d$.

\begin{theorem}[{\cite[Thm.~1.2, Cor.~5.15]{KinjoKoseki}}]
For any $r$ and any two $d,d'$ there is a canonical isomorphism
$$	h_*(\mathcal{BPS}^C({r,d})) \cong h_*(\mathcal{BPS}^C({r,d'})) .$$
Passing to global sections, this yields isomorphisms of graded vector spaces
\begin{equation}\label{eq:chi_indepBPS}
\widehat{\g}^C(r,d) \cong \widehat{\g}^C(r,d'), \quad {\g}^C(r,d) \cong {\g}^C(r,d') 
\end{equation}
which are compatible with the perverse filtration.
\end{theorem}

\medskip\subsection{Stacky P=C} The P=C equality admits an extension to the case of an arbitrary pair $(r,d)$ due to R. Picovschi. 
It is best understood in terms of the (affinized) BPS Lie algebra. Recall the ring
$$S_C=\Q[c_k(\alpha)\;;\; k \geqslant 1, \alpha] =Q[\ch_k(\alpha)\;;\; k \geqslant 1, \alpha] ,$$
where $\alpha$ runs over a basis of $\H^*(C,\Q)$.
We equip $S_C$ with the \textit{Chern} grading defined by
$$\deg(\ch_k(\alpha))=k \quad \forall \; k \geqslant 1, \alpha \in \H^*(C,\Q).$$
Notice that $\deg(\ch_k(\alpha))$ coincides with the cohomological degree (up to a factor of 2) if and only if $\alpha \in \Q [\text{pt}]$.
Let us write 
$$S_C= \bigoplus_{n \geq 0}S_{C,n}, \qquad S_{C,\leqslant n}=\bigoplus_{\ell \leqslant n} S_{C,n} .$$
For every $(r,d)$ we define an automorphism of $S_C$ given on generators by
$$\varphi_{r,d}: \ch_n(\alpha)\mapsto\ch_n(\alpha) + \frac{dn}{r}\ch_{n-1}([\text{pt}]\cup \alpha), \qquad (n \geqslant 1, \;\alpha) .$$
Note that $\varphi_{r,n} \circ \varphi_{r,n'}=\varphi_{r,n+n'}$, $\varphi_{r,0}=\text{Id}$, $\varphi_{r,d}(S_{C,\leqslant n})=S_{C,\leqslant n}$. Note also that $\varphi_{r,d}(\ch_k(\alpha))=\ch_k(\alpha)$ unless
$|\alpha|=0$.
\begin{theorem}[{\cite{Raphael}}]\label{stackyP=C}
There exists a splitting of the perverse filtration
$$\P_k\,\widehat{\g}^C(r,d)=\bigoplus_{\ell=0}^{k} L_\ell\,\widehat{\g}^C(r,d)$$
such that 
$$\varphi_{r,d}(\overline{c}_s(\eta)) \bullet  L_\ell\, \widehat{\g}_k^C(r,d) \subseteq L_{\ell+s}\, \widehat{\g}_{k-2s+2-|\eta|}^C(r,d)
, \qquad $$
$$
u \bullet  L_\ell \,\widehat{\g}_k^C(r,d) \subseteq L_{\ell+2} \,\widehat{\g}_{k-2}^C(r,d)$$
for any $s \geqslant 1$ and $(s,\eta) \neq (1,[\text{pt}])$.
\end{theorem}

As a consequence, we have
\begin{equation}\label{eq:tautperversepureHiggs}
\bigoplus_\ell L_{N-\ell}\, \widehat{\g}_{2\ell}^C(r,d)=\Q[ \overline{c}_2([\text{pt}]), \overline{c}_3([\text{pt}]), \ldots ] \bullet [\mathfrak{M}^{\ss,\rd}(r,d)]
\end{equation}
and for any sequence $(s_1,\ldots, s_p)$ we have
\begin{equation}\label{eq:tautperverspureHiggs2}
\prod_i \overline{c}_{s_i}([\text{pt}]) \cap [\mathfrak{M}^{\ss,\rd}(r,d)] \in L_{\sum_i s_i}\widehat{\g}^C(r,d).
\end{equation}

\begin{remark}
For the purpose of this paper, only the pure part of Theorem~\ref{stackyP=C}, i.e. \eqref{eq:tautperversepureHiggs} and \eqref{eq:tautperverspureHiggs2} are needed. This statement can be given a simpler proof, based on \cite{HMMS} and \cite{DHKSV}. For the reader's comfort, we provide a proof of these in Appendix~\ref{app:C}. 
We are greatful to Jiahao Bi for pointing this out.
\end{remark}

%=====================================================================
\medskip\section{Restriction around the trivial Higgs bundle.}
%=====================================================================

\medskip

We keep the notation of the previous section: $C$ is a smooth connected projective complex curve of genus $g$,
$\mathfrak{M}^\ss(r,0)$ and $M^D(r,0)$ are the derived moduli stack and coarse moduli space of semistable Higgs bundles on $C$ of rank $r$ and degree $0$,
and $p: \mathfrak{M}^{\ss,\cl}(r,0) \to M^D(r,0)$ is the good moduli space.

\medskip\subsection{Davison's neighborhood theorem}\label{sec:restriction_map} The Higgs bundle $(\mathcal{O}^r,\varphi=0)$ is polystable. 
Thus it corresponds to a point of $M^D(r,0)$. Let $\calO^r$ denote this point. We have 
$$\text{Ext}^0(\calO,\calO) \cong \text{Ext}^2(\calO,\calO)^\vee \cong \C
,\quad
\text{Ext}^1(\calO,\calO) \cong T^*\H^1(C,\calO^r)=\C^{2g}$$ 
The associated $\text{Ext}$-quiver is the $S_g$ quiver. 
The relation between the COHA and the BPS Lie algebra $\mathbf{A}^{C,(0)}$ and
$\mathfrak{g}^{C,(0)}$ attached to $\mathfrak{M}^\ss(r,0)$, and the COHA  and BPS Lie algebra 
$\mathbf{A}^{S_g}$ and $\mathfrak{g}^{S_g}$ attached to $S_g$ is given by the following theorem due to Davison.

\begin{theorem}[{\cite[Thm.~5.11]{D24}}] 
 There are open analytic neighborhoods $U$ of $\calO^r$ in $M^D(r,0)$ and $U'$ of $0$ in 
 $\text{Rep}_r\Pi_{S_g}$ together with analytic isomorphisms $\phi: U \cong U'$ mapping $\calO^r$ to $0$ and 
 $\boldsymbol{\phi}: \mathfrak{U} =p^{-1}(U) \cong p^{-1}(U')=:\mathfrak{U}'$ making the following diagram commute
$$
\begin{tikzcd}
\mathfrak{M}^{\ss,\cl}(r,0) \arrow[d,"p"] & \mathfrak{U}\arrow[d,"p"] \arrow[l,"j_{\mathfrak{U}}"'] \arrow[r,"\boldsymbol{\phi}"]
& \mathfrak{U}' \arrow[r,"j_{\mathfrak{U}'}"] \arrow[d,"p"] & \Rep^\cl_r\Pi_{S_g} \arrow[d,"p"]\\ M^D(r,0) & U \arrow[r,"\phi"] \arrow[l,"j_U"'] 
& U' \arrow[r,"j_{U'}"] & \text{Rep}_r\Pi_{S_g}. 
\end{tikzcd}
$$
Here the leftmost and rightmost horizontal arrows are open embeddings.
\end{theorem}

Note that, since $\mathfrak{Rep}_r\Pi_{S_g}$ is a conical stack with apex $0$, 
the restriction to any star-shaped open neighborhood $\mathfrak{U}'$ 
of $0$ induces an isomorphism in Borel-Moore homology 
$$\H_*(\mathfrak{U}',\Q) \cong \H_*(\mathfrak{Rep}_r\Pi_{S_g},\Q) .$$
Composing with $\phi_*$, we obtain a restriction morphism, which is independent of the choices of $U,$ $U'$:
$$\text{res}_U: \H_*(\mathfrak{M}^\ss(r,0),\Q) \to \H_*(\Rep_r\Pi_{S_g},\Q).$$
Considering all values of $r\geqslant 0$, we can form a system of open neighborhoods $U_r$ in $M^D(r,0)$ such that
$$\oplus_{r,d}^{-1}(U_{r+s}) \subset U_r \times U_s,\quad \forall r,s \geqslant 0.$$
The same thus holds for $\{U'_r\}_r$ in $\text{Rep}_r\Pi_{S_g}$.

\begin{corollary}
The collection of restriction maps $\text{res}_{U_r}$ induces an algebra morphism $\mathbf{A}_C^{(0)} \to \mathbf{A}_{\Pi_{S_g}}$ 
which restricts to Lie algebra morphisms $\g^{C,(0)} \to \g^{\Pi_{S_g}}$ and $\widehat{\g}^{C,(0)} \to \widehat{\g}^{\Pi_{S_g}}$.
\end{corollary}
\begin{proof}
The COHA multiplication $\mathbf{A}_C(r,0) \otimes \mathbf{A}_C(s,0) \to \mathbf{A}_C(r+s,0)$ 
is obtained by taking derived global sections of a morphism of complexes over $M^D(r+s,0)$
$$\oplus_* \left(p_*(\DQ_{\mathfrak{M}^\ss(r,0)} \otimes \mathbb{L}^{(g-1)r^2})\boxtimes p_*(\DQ_{\mathfrak{M}^\ss(s,0)} \otimes \mathbb{L}^{(g-1)s^2}) \right) \to p_*(\DQ_{\mathfrak{M}^\ss(r+s,0)} \otimes \mathbb{L}^{(g-1)(r+s)^2}) $$
and likewise for $\mathbf{A}_{\Pi_{S_g}}$. By open base change and further restricting to $\mathfrak{U}_{r+s}$ we get a similar morphism
$$\oplus_* \left(p_*(\DQ_{\mathfrak{U}_r} \otimes \mathbb{L}^{(g-1)r^2})\boxtimes p_*(\DQ_{\mathfrak{U}_s} \otimes \mathbb{L}^{(g-1)s^2}) \right) \to p_*(\DQ_{\mathfrak{U}^{r+s}} \otimes \mathbb{L}^{(g-1)(r+s)^2})$$
and likewise on the quiver side. The restriction map is thus an algebra morphism on the Hitchin side and an isomorphism on the quiver side. 
This proves the first statement. 
The second statement follows from the fact that the construction of the BPS sheaf is local over the coarse moduli space.
\end{proof}

\medskip\subsection{Restriction of tautological classes} 
Davison's neighborhood isomorphism is compatible with the reduced virtual fundamental classes in the sense that
\begin{equation}\label{eq:resredvirfundclass}
\text{res}_U([\mathfrak{M}^{\ss,\rd}(r,0)])=[\Rep^\rd_r\Pi_{S_g}] .
\end{equation}
As we could not locate this statement in the literature, we provide an argument when $g>1$,
which is enough for our purpose. In that case, we claim that $[\mathfrak{M}^{\ss,\rd}(r,0)]$ and $[\Rep^\rd_r\Pi_{S_g}]$ are both equal to the 
corresponding classical fundamental classes, and these are obviously compatible with open restriction. Indeed, by 
\cite[Prop.~8.1]{MMSV} the stack $\mathfrak{M}^\ss(r,0)$ is irreducible and of classical dimension $2N-1$, which is the same as the 
virtual dimension of $\mathfrak{M}^{\ss,\rd}(r,0)$. Likewise, by \cite[Cor.~1.4.]{CB}, the stack $\Rep^\cl_r\Pi_{S_g}=\mu_r^{-1}(0)/GL_r$ is of 
dimension $2N-1$, which is the same as the virtual dimension of $\Rep^\rd_r\Pi_{S_g}$. In addition, both stacks have dense open 
substacks which are smooth: it suffices to take the respective substacks parametrizing simple sheaves or representations. 
Hence, by Corollary~\ref{cor:vir_generically_smooth}, it follows that
$$[\mathfrak{M}^{\ss,\rd}(r,0)]=[\mathfrak{M}^{\ss,\cl}(r,0)], \quad [(\Rep^\cl_r\Pi_{S_g}]=[\Rep^\rd_r\Pi_{S_g}]$$
as wanted.

We next examine the behavior of $\text{res}_U$ with respect to tautological classes. 
By construction, the tautological classes on $\Rep_r\Pi_{S_g}$ are all pulled back via the morphism
$\Rep_r\Pi_{S_g} \to BGL_r(\mathbb C)$.
We denote by 
$$\ch_k=\ch_k(V) \in \H^{2k}(\Rep_r\Pi_{S_g},\Q)$$ 
the $k$th Chern character of the tautological bundle $V$. 
On the Higgs moduli side, we fix a point $x\in C$ and let 
$\calV_x$ be the restriction of the tautological bundle $\calV$ to $\mathfrak{M}^\ss(r,0) \times \{x\}$.
It induces a morphism $$\epsilon_x\colon\mathfrak{M}^\ss(r,0)\to BGL_r(\mathbb C) .$$
We consider the closed immersions
$$i_{\calO^r}: BGL_r(\mathbb{C}) =B\text{Aut}(\calO^r) \to \mathfrak{U}, \quad i_{0}: BGL_r(\mathbb{C}) =B\text{Aut}(0) \to \mathfrak{U}' $$ 
 Finally, let us denote by $$u: \mathfrak{U}\to \mathfrak{M}^\ss(r,0), \quad u': \mathfrak{U}' \to \Rep_r\Pi_{S_g}$$ 
the open analytic embeddings.

\begin{lemma}\label{lem:ch-eta}
The following hold: 
\begin{enumerate}[label=$\mathrm{(\alph*)}$,leftmargin=8mm,itemsep=1mm]
	\item The maps $i_{0}^*:\H^*(\mathfrak{U}',\Q) \to \H^*_{GL_r}$ and $(u')^*: \H^*(\Rep_r\Pi_{S_g},\Q) \to \H^*(\mathfrak{U}',\Q)$ are 
	algebra isomorphisms.
	\item The maps $\epsilon_x^*: \H^*_{GL_r} \to \H^*(\mathfrak{U},\Q)$ and $i_{\calO^r}^*:\H^*(\mathfrak{U},\Q) \to \H^*_{GL_r}$ are 
	inverse algebra isomorphisms.
\end{enumerate}
Composing these maps yields an algebra morphism $\text{res}_U: \H^*(\mathfrak{M}^\ss(r,0),\Q) \to \H^*(\Rep_r\Pi_{S_g},\Q)$.
\begin{enumerate}[label=$\mathrm{(\alph*)}$,leftmargin=8mm,itemsep=1mm]
	\item[$\mathrm{(c)}$] We have $\text{res}_U(\ch_k([\text{pt}]))=\ch_k$ for any $k \geqslant 1$.
	\item[$\mathrm{(d)}$] The map $\text{res}_U$ restricts to a graded algebra isomorphism $\Q[\ch_1([\text{pt}]), \ch_2([\text{pt}]), \ldots] \cong \H^*(\Rep_r\Pi_{S_g},\Q)$ .
\end{enumerate}
\end{lemma}

\begin{proof} Part (a) comes from the fact that $\Rep_r\Pi_{S_g}$ is conical, $0$ is its apex and $\mathfrak{U}'$ is a star-shaped open neighborhood. Next, the composition
	\begin{equation}\label{eq:cohU}
		\begin{tikzcd}
			BGL_r(\mathbb{C}) =B\text{Aut}(\calO^r) \arrow[r,"i_{\calO^r}"] & \mathfrak{U} \arrow[r,"\epsilon_x"] & BGL_r(\mathbb{C})
		\end{tikzcd}
	\end{equation}
is the identity morphism, hence $i_{\calO^r}^*\circ \epsilon_x^* =\Id$. 
By (a), the map $i_{\calO^r}^*$ is an isomorphism and (b) follows. As $[\text{pt}] \in \H^2(C,\Q)$ is Poincar\'e dual to a point, we have 
$$\ch_k([\text{pt}])=\int_C \ch_k(\calV) \cup p_2^*([\text{pt}])=\iota_x^*\ch_k(\calV)=\ch_k(\calV_x)=\epsilon_x^*(\ch_k),$$
where $\iota_x: \{x\} \to C$ is the inclusion. Statement (c) easily follows. Finally (d) is a consequence of (c) and of the fact that $\H^*(\Rep_r\Pi_{S_g},\Q) \cong \H^*_{GL_r}$. 
\end{proof}

Note that the restriction maps $\text{res}_U$ in Borel-Moore homology and in cohomology are compatible with the $\bullet$ action, i.e., 
for any $c \in \H^*(\mathfrak{M}^{\ss}(r,0),\Q)$ and $a \in \H_*(\mathfrak{M}^{\ss}(r,0),\Q)$, we have
\begin{equation}
\text{res}_U(c \bullet a)=\text{res}_U(c) \bullet \text{res}_U(a) .
\end{equation}

\medskip\section{Proof of the Hausel-Rodriguez-Villegas purity conjecture.}

\medskip

We now have all the ingredients needed to prove Theorem~\ref{thm:purityconj}. 
Fix a coprime pair $(r,d)$. Recall that $N=(g-1)r^2+1$. The P=W theorem yields
$$\sum_k \dim\;\text{W}_{2k}\,\H^{2k}(M^B(r,d),\Q)q^k= \sum_k \dim\;\gr^\P_k\,\H^{2k}(M^D(r,d),\Q)q^k$$
and by \eqref{eq:lowerbound} we have
\begin{equation}\label{eq:finalproof1}
\sum_k \dim\;\gr^\P_k\,\H^{2k}(M^D(r,d),\Q)=A_{S_g,r}(1) .
\end{equation}
On the other hand, by Poincar\'e duality, see \eqref{eq:PoincaréBPS}, \eqref{eq:isomgradedpureperverse} 
and the $\chi$-independence \eqref{eq:chi_indepBPS} we have
\begin{equation}\label{eq:finalproof2}
	\begin{split}
	\sum_k \dim\;\gr^\P_k\,\H^{2k}(M^D(r,d),\Q)q^k &= \sum_\ell \dim\;\gr^\P_{N-\ell}\, \g^C_{2\ell}(r,d)q^{N-\ell}\\
	&=\sum_\ell \dim\;\gr^\P_{N-\ell}\, \widehat{\g}^C_{2\ell}(r,d)q^{N-\ell} \\
	&=\sum_\ell \dim\;\gr^\P_{N-\ell}\, \widehat{\g}^C_{2\ell}(r,0)q^{N-\ell} \\
	&=\sum_\ell \dim\;L_{N-\ell}\, \widehat{\g}^C_{2\ell}(r,0)q^{N-\ell} .
	\end{split}
\end{equation}
Next, by the stacky P=C theorem, see \eqref{eq:tautperversepureHiggs}, we have
\begin{equation}\label{eq:finalproof3}
	\dim\;L_{N-\ell}\,\widehat{\g}^C_{2\ell}(r,0)=\dim\;\overline{S}_{C,N-\ell} \bullet [\mathfrak{M}^{\ss,\rd}(r,0)]
\end{equation}
where 
$$\overline{S}_C=\bigoplus_k \overline{S}_{C,k}=\Q[ \overline{c}_2([\text{pt}]), \overline{c}_3([\text{pt}]), \ldots, \overline{c}_r([\text{pt}])] \subset S_C $$
and $\deg(\overline{c}_k([\text{pt}]))=k$ for $k=2,\ldots, r$. 
Applying the restriction map $\text{res}_U$ of \S\ref{sec:restriction_map} and using \eqref{eq:resredvirfundclass} together with Lemma~\ref{lem:ch-eta} (c), we deduce that for any $\ell$ we have
\begin{equation}\label{eq:finalproof4}
\dim\;\overline{S}_{C,N-\ell} \bullet [\mathfrak{M}^{\ss,\rd}(r,0)] \geqslant \dim\;\overline{S}_{S_g,N-\ell} \bullet [\Rep^\rd_r\Pi_{S_g}] ,
\end{equation}
where
$$\overline{S}_{S_g}=\bigoplus_k \overline{S}_{S_g,k}=\Q[\overline{c}_2, \overline{c}_3, \ldots, \overline{c}_r] \subset S_{S_g} $$
and $\deg(\overline{c}_k)=k$ for $k=2,\ldots, r$. By the quiver version of the Markman Theorem \ref{thm:Markman2} we have
$$\overline{S}_{S_g}[c_1] \bullet [\Rep^\rd_r\Pi_{S_g}]=\widehat{\g}^{\Pi_{S_g}}(r).$$
Because $\Q[c_1]$ acts freely on $\widehat{\g}^{S_g}(r)$ and 
$\widehat{\g}^{S_g}(r)/c_1\bullet \widehat{\g}^{S_g}(r) \cong {\g}^{\Pi_{S_g}}(r)$, see \eqref{eq:gfromghat}, we deduce that
$$\dim\;\overline{S}_{S_g,N-\ell}\bullet [\Rep^\rd_r\Pi_{S_g}]\geqslant\dim\;\g^{\Pi_{S_g}}_{2\ell}(r)$$
and hence from \eqref{eq:dimBPSLiealg} that
\begin{equation}\label{eq:finalproof5}
\sum_\ell \dim\;\overline{S}_{S_g,N-\ell} \bullet [\Rep^\rd_r\Pi_{S_g}] q^{N-\ell}\geqslant q^NA_{S_g,r}(q^{-1}) ,
\end{equation}
where $\geqslant$ means coefficient-wise inequality.
Combining \eqref{eq:finalproof2}, \eqref{eq:finalproof3}, \eqref{eq:finalproof4} and \eqref{eq:finalproof5} we see that
$$\sum_k \dim\;\gr^\P_k \,\H^{2k}(M^D(r,d),\Q)q^k \geqslant q^NA_{S_g,r}(q^{-1}) .$$
Because of \eqref{eq:finalproof1}, it must be an equality: indeed, two polynomials $P,P' \in \mathbb{N}[q]$ for which $P \geqslant P'$ and $P(1)=P'(1)$ must be equal. This finishes the proof of Theorem~\ref{thm:purityconj} .

\bigskip

\centerline{\textbf{AI disclosure}}

\smallskip

The paper was fully written by us and all the proofs and ideas are ours. AI (Claude and ChatGPT) 
was used for proofreading purposes and to improve the presentation. 

\medskip

\centerline{\textbf{Acknowledgements}}

\smallskip

We are grateful to J. Bi, B. Davison, L. Hennecart, T. Kinjo, R. Picovschi and N. Seroux for fruitful discussions and correspondence. The quiver version of Markman's theorem was first announced by Hennecart in the workshop 'Algebraic Geometry Days' 2026. The proofs were worked out independently and communicated between us and the authors of \cite{DH} before submission to the arXiv.

\medskip

\appendix

%=====================================================================
\medskip\section{Complements on the coproduct on $\mathbf{A}_{\Pi_Q}$}\label{app:A}
%=====================================================================

The aim of this Appendix is to collect and adapt from \cite[\S 2, \S 6]{DHKSV} to the present quiver context some known results and constructions concerning the BPS Lie algebra, its coproduct and the action of tautological classes.

\medskip\subsection{Transfer of the construction of the coproduct \cite{DHKSV} to quivers}\label{sec:transfer}

We let $Q=(Q_0,Q_1)$ be an arbitrary finite quiver and use the notations from \S~2 regarding the stacks of representations of $Q$ and $\Pi_Q$. We keep the notations from the main body of the article.

\iffalse%%%%%%%%
\smallskip\subsubsection{The $(-1)$-shifted side}
The $(-1)$-shifted cotangent stack 
$$\tcalM_d^Q=T^*[-1]\calM_d^Q$$ is the derived stack of $d$-dimensional modules over the Jacobi algebra of the tripled quiver $\widetilde Q$
 (=the double $\overline Q$ with a loop $\omega_i$ added at each vertex) 
 with potential $W=\omega\sum_a[x_a,x_a^*]$. That is, of pairs $(M,\omega)$ 
 with $M$ a $\Pi_Q$-module and $\omega\in\operatorname{End}(M)$. 
 It carries the standard $(-1)$-shifted symplectic structure and orientation of \cite[\S2.1]{DHKSV}.
 The DT sheaf $\varphi_{\tcalM^Q}$ is the sheaf of vanishing cycles of $\tr W$ following \cite{BBDJS15}. 
 Dimensional reduction \cite[\S2.2.3]{DHKSV} identifies
\[
\H^*(\tcalM_d^Q,\varphi_{\tcalM_d^Q})\cong\H^\BM_{-*+\vdim\calM_d^Q}(\calM_d^Q,\Q)=\bfH^Q_d,
\]
compatibly with the COHA products by [BenDimRed].

\fi%%%%%%%%%%%%

\smallskip\subsubsection{Hypotheses of \cite[\S2.4.1]{DHKSV}}
The following hold, for any dimension vector $d \in \N Q_0$.
\begin{enumerate}[label=(\alph*),leftmargin=8mm,itemsep=2mm]
\item $\Rep_d\Pi_Q$ is a quasi-smooth quasi-compact Artin $1$-stack, being the quotient of a derived affine scheme by $GL_d$.
\item $\bigsqcup_d \Rep_d\Pi_Q$ is closed under direct sums and summands, and $\Rep_0\Pi_Q=\{0\}$.
\item $\Rep_d\Pi_Q$ has affine diagonal (it is a global quotient of an affine scheme by an affine group).
\item $\Rep_d\Pi_Q$ is $\Theta$-reductive: the stack $\mathfrak{Filt}_{d,d'}\Pi_Q$ is a closed substack of a bundle of quiver flag varieties, hence the morphism $\mathfrak{Filt}_{d,d'}\Pi_Q\to\Rep_{d+d'}\Pi_Q$ is proper .
\item Each $\Rep_d\Pi_Q$ is quasi-compact. The graded points are quasi-compact for the same reason.

\item Each $\Rep_d\Pi_Q$ has a quasi-compact good moduli space $p: \Rep^\cl_d\Pi_Q \to Rep_d\Pi_Q$.
\end{enumerate}
 
The $2$CY property of $\Pi_Q$-mod, which fails for Dynkin $Q$, is used in \cite{DHKSV} only to produce the 
$0$-shifted symplectic structure on $\Rep_d\Pi_Q$ 
and the orientation data on $\tRep_d\Pi_Q$. Here both come for free from the cotangent description. 
Note that, unlike \cite{DHKSV}, each $\Rep_d\Pi_Q$ being quasi-compact, no pro-objects are needed. 
Further, we do not consider categories of semistable objects, or, rather, we only consider the zero stability condition, i.e.,
for which every object is semistable.

The above hypotheses (a)-(e) imply the existence of a COHA structure on $\mathbf{A}_{\Pi_Q}$. Condition (f) implies the existence of a local COHA, BPS sheaves $\mathcal{BPS}_d \in \text{MHM}(Rep_d\Pi_Q)$ and a BPS (resp. affinized BPS) Lie algebra $\g^{\Pi_Q}$ 
(resp. $\widehat{\g}^{\Pi_Q}$) satisfying the absolute integrality (PBW) theorem. More precisely, the inclusion 
$$\mathcal{BPS}^{\Pi_Q}(d) \to p_*(\DQ_{\mathfrak{M}_d\Pi_Q})$$
 is defined by dimensional reduction from the inclusion of the lowest perverse cohomology sheaf $${}^\mathrm{p}\mathcal{H}^1\left((\tilde p_d)_*\varphi_{\tRep_d\Pi_Q}\right)$$ in $(\tilde p_d)_*\varphi_{\tRep_d\Pi_Q}$, the BPS Lie algebra is then defined as
$$\g^{\Pi_Q}(d)= \H^{-*}\Big(Rep_d\Pi_Q, \mathcal{BPS}^{\Pi_Q}(d)\Big) .$$
To define the affinized BPS Lie algebra, we consider the determinant line bundle
$$\mathcal{L}_d=\bigotimes_{i \in Q_0} \det V_i$$
where $V_i \in Coh(\Rep_d\Pi_Q)$ is the tautological rank $d_i$ vector bundle corresponding to the $i$th vertex. Under the standard $\BGM$-action, the weight of $\det V$ is non-zero, so that $\det V$ is an equivariant parameter in the sense of \cite[\S2.5.2]{DHKSV}. Slightly abusing notation, we put 
$$u=c_1(\mathcal{L}_d) \in \H^2_{GL_d}$$
which we consider as an operator via the cap product on $\bfA_{\Pi_Q}$. We now define $\widehat{\g}^{\Pi_Q}(d)$ as the image of the subspace $\g^{\Pi_Q}(d) \otimes \H^*_{\GM}$ under the map
$$\bfA_{\Pi_Q}(d) \otimes \H^*_{\GM} \to \bfA_{\Pi_Q}(d), \quad x \otimes z^\ell \mapsto u^\ell \bullet x .$$

The following microlocal characterization of the affinized BPS Lie algebra is proved as in \cite[Cor.~2.25]{DHKSV}:

\begin{lemma} We have
	\begin{equation}\label{eq:g-BPS}
		\widehat{\g}^{\Pi_Q}(d)
		=\Im\Big(\H^{\Lambda_\snilp}_*(\Rep_d\Pi_Q,\Q)\to\mathbf{A}_{\Pi_Q}(d)\Big).
	\end{equation}
	\qed
\end{lemma}

\smallskip\subsubsection{Coproduct and Milnor-Moore theorem}
For an analytic open $U\subset\C$, let 
$$\tRep_U\Pi_Q\subset\bigcap_d\tRep_d\Pi_Q$$ be the open substack of pairs $(M,\theta)$ with all eigenvalues of $\theta$ in $U$. 
The openness proof of \cite[\S6.2]{DHKSV} applies verbatim. 
From this, the construction of \cite[\S6.2]{DHKSV} produces a coassociative, cocommutative coproduct 
$\Delta$ on $\mathbf{A}_{\Pi_Q}$ which is compatible with the product, making $\bfA_{\Pi_Q}$ into a connected, $\N Q_0\times\Z$-graded
 cocommutative bialgebra with finite-dimensional graded pieces. The Euler form $(-,-)$ is symmetric, so the sign-twisted braiding of \cite[\S 6.1]{DHKSV} is defined.
By the Milnor-Moore theorem, we have
\begin{equation}\label{eq:MM}
\bfA_{\Pi_Q}\cong\U(\bfA^{\text{prim}}_{\Pi_Q}),\quad\bfA^{\text{prim}}_{\Pi_Q}=\{x\in\bfA_{\Pi_Q}\;;\;\Delta(x)=x\otimes1+1\otimes x\}.
\end{equation}
The proof of \cite[Lem.~6.12]{DHKSV} applies verbatim to give the following 
identification of $\widehat{\g}^{\Pi_Q}$ with the Lie algebra of primitive elements.

\begin{lemma} We have $\widehat{\g}^{\Pi_Q}(d)=\mathbf{A}_{\Pi_Q}^\mathrm{prim}(d)$ for any $d \in \N Q_0$.
\end{lemma}

%%%%%%%%%%%%%%%%%%%%

\medskip\subsection{Sketch of the proof of Lemma~\ref{lem:hopf-module}} For simplicity, we drop the mention of the dimension vector $d$ in 
the indices. 
We begin with (a). Recall the induction diagram \eqref{diag:defprodcoha1}  defining the COHA multiplication. By additivity of the Chern character
$$\mathrm{ev}^*(c)=\mathrm{gr}^*( \Delta c),\quad c\in S_Q.$$
The projection formula then gives
$$
c\bullet\mathrm{ev}_*\mathrm{gr}^!(x\otimes y)=\mathrm{ev}_*\big(\mathrm{gr}^*(\Delta( c))\cap\mathrm{gr}^!(x\otimes y)\big)=\mathrm{ev}_*\mathrm{gr}^!\big(\Delta(c)\bullet(x\otimes y)\big).
$$
Now, we turn to (b). For $U \subset \mathbb{A}^1$ an analytic open subset, put
$$(\tRep\Pi_Q)_U=\{(x,x^*,\theta)\in \tRep\Pi_Q\,;\,\Spec(\theta) \subset U\}
,\quad 
\DT(U)=\H^*((\tRep\Pi_Q)_U,\varphi_{\Rep\Pi_Q}|_{(\tRep\Pi_Q)_U}).$$ 
The coproduct is constructed via a factorization process and is given by the restriction $\DT(U)\to\DT(U_1\sqcup U_2)$ of a the so-called Donaldson-Thomas sheaf and by the 
identification 
\begin{equation}\label{eq:product factorization}
	(\tRep\Pi_Q)_{U_1\sqcup U_2}\cong (\tRep\Pi_Q)_{U_1}\times(\tRep\Pi_Q)_{U_2} .
\end{equation}
Under the dimensional reduction, cap product with $c$ on $\bfA_{\Pi_Q}$ corresponds to the action of 
$$\kappa^*(c)\in\H^*(\tRep\Pi_Q,\Q)$$ on $\H^*(\tRep\Pi_Q,\varphi_{\tRep\Pi_Q})$.
Restriction commutes with the action of cohomology classes, and 
under the identification \eqref{eq:product factorization}, $\kappa^*c$ restricts to $\sum \kappa^*c_{(1)}\boxtimes \kappa^*c_{(2)}$, where $\Delta(c)=\sum c_{(1)} \otimes c_{(2)}$. We have  
again used additivity of the Chern character, 
this time applied to $V=V_{U_1}\oplus V_{U_2}$.

Finally, we deal with (c). If $x$ is primitive, then by (b) we have 
$$\Delta(c\bullet x)=\sum (c_{(1)}\bullet x)\otimes(c_{(2)}\bullet1)+\sum(c_{(1)}\bullet1)\otimes(c_{(2)}\bullet x).$$
Since $\bfA_{\Pi_Q,0}=\Q$ and $S_Q$ acts on it through the counit, the rhs is equal to 
$$(c\bullet x)\otimes1+1\otimes(c\bullet x)$$ as wanted. \qed

%%%%%%%%%%%%%%%%

\medskip\section{Quasi-smooth stacks with smooth truncation}

\medskip\subsection{Conventions}

We work over a field $k$ of characteristic 0. Derived stacks are derived Artin stacks locally almost
of finite presentation over $k$, modelled on commutative simplicial algebras or dg-algebras.
Complexes are graded {cohomologically}, so $h^{j}(-)=\pi_{-j}(-)$. A
complex $F$ is {$n$-connective} if $h^j(F)=0$ for $j>-n$. It is connective if it is 0-connective.
Let $A$ be a connective simplicial commutative ring. An $A$-module $N$ is {discrete} if its homotopy is
concentrated in degree $0$.
Likewise, a quasi-coherent sheaf on a derived stack $Y$ is {discrete} if its
cohomology is concentrated in  degree $0$.
If $N$ is a discrete $A$-module, the $A$-action factors through $\pi_0A$: the
higher homotopy groups $\pi_{n}A$, $n\geqslant 1$, act through maps
$\pi_nA\otimes\pi_0N\to\pi_nN=0$. The functor $N\mapsto\pi_0N$ gives an
equivalence
\[
\{\text{discrete }A\text{-modules}\}\;\cong\;\{\text{ordinary }\pi_0A\text{-modules}\}
\]
Similarly, discrete
quasi-coherent sheaves on $Y$ are the same as ordinary quasi-coherent sheaves on
the classical truncation $t_0Y$.

A perfect complex $E$ on a derived stack $Y$ has {Tor-amplitude in $[a,b]$}
if for every discrete quasi-coherent sheaf $N$ on $t_0Y$ (or, smooth locally,
every discrete module $N$ over $\pi_0 A$ for an affine chart $\Spec A\to Y$)
the complex $E\otimes^{\L} N$ has cohomology concentrated in degrees $[a,b]$.
If $E$ has Tor-amplitude in $[a,b]$, then $h^j(E)=0$ for $j>b$. If moreover $Y$
is classical, then (taking $N=\calO_Y$) the cohomology sheaves of $E$ are
concentrated in $[a,b]$. This fails for derived $Y$: for $Y=\Spec A$ with
$A=k[\varepsilon]$ and $\varepsilon$ of homological degree $1$, the complex $\calO_Y$ has
Tor-amplitude in $[0,0]$ but $h^{-1}(\calO_Y)=\pi_1A=k\neq0$.
Tor-amplitude in $[0,0]$ means
$\Tor_j(E,N)=0$ for $j>0$, i.e., $E$ is flat, hence, since $E$ is perfect it is a vector
bundle in degree $0$. Tor-amplitude in $[-1,0]$ means that $E$ is locally
quasi-isomorphic to a two-term complex of vector bundles $E^{-1}\to E^{0}$.
For a perfect complex $E$ on a connected $Y$ let $\chi(E)\in\mathbb{Z}$ be its
rank (=Euler characteristic).

A derived stack $X$ is {quasi-smooth} if its cotangent complex $\L_X$ is
perfect of Tor-amplitude in $[-1,1]$. Its {virtual dimension} is the
locally constant function $\vd X=\chi(\L_X)$.
Let $i\colon X^\cl\hookrightarrow X$ be the inclusion of the classical truncation.
The complex $E=i^*\L_X$ is perfect of Tor-amplitude in $[-1,1]$, and we let
\[
\phi\colon E\longrightarrow\L_{X^\cl}
\]
be the natural map. It is the (perfect) obstruction theory induced on $X^\cl$
\cite{STV}.

\begin{theorem}\label{thm:mainB}
Let $X$ be a quasi-smooth derived stack. Suppose that $M=X^\cl$ is smooth over
$k$ with $\dim M=\vd X$ on every connected component. Then
\begin{enumerate}[label=$\mathrm{(\alph*)}$,leftmargin=8mm,itemsep=1mm]
  \item $\phi\colon E\to\L_M$ is an equivalence,
  \item $i\colon M\to X$ is an equivalence, 
  \item the virtual fundamental class of $X$ is the fundamental class, i.e.,
        $[X]^{\vir}=[M]$.
\end{enumerate}
\end{theorem}

Throughout we assume $M$ is connected, so that $\dim M$ and $\vd X$ are integers.

\medskip\subsection{Lemmas}

\begin{lemma}\label{lem:conn}
Let $f\colon A\to B$ be a map of connective commutative dg-algebras over $k$.
Suppose the mapping cone
$\cofib(f)$ is $m$-connective for some $m\geqslant 1$. Then $\L_{B/A}$ is
$m$-connective, and the natural map
\[
B\otimes_A\cofib(f)\longrightarrow \L_{B/A}
\]
has $(m+1)$-connective mapping cone. In particular
$$\pi_m\L_{B/A}\cong \pi_m\big(B\otimes_A\cofib(f)\big)\cong \pi_0B\otimes_{\pi_0A}\pi_m\cofib(f).$$
\end{lemma}

\begin{proof}
This is \cite[Thm.~7.4.3.1]{HA} and, in the simplicial setting, \cite[\S2.2.2]{TV08}.
\end{proof}

\begin{lemma}\label{lem:Lt0}
For any derived stack $X$, the relative cotangent complex $\L_{X^\cl/X}$ is $2$-connective. Consequently
$\phi\colon E\to\L_{X^\cl}$ induces isomorphisms on $h^j$ for $j\geqslant 0$ and
a surjection on $h^{-1}$.
\end{lemma}

\begin{proof}
The question is smooth local on $X$ because for a smooth chart $U\to X$ we have
$U^\cl=U\times_X X^\cl$ and $\L_{U^\cl/U}$ is the pullback of $\L_{X^\cl/X}$.
So we may assume
$$X=\Spec A,\quad X^\cl=\Spec \pi_0A,\quad f\colon A\to\pi_0A.$$
The mapping cone of $f$ is
$\cofib(f)=\tau_{\geqslant 1}A[1]$. It is $2$-connective.
Lemma~\ref{lem:conn} with $m=2$ gives the first claim. The second follows from
the long exact sequence of the sequence
$E\to\L_{X^\cl}\to\L_{X^\cl/X}$.
\end{proof}

\begin{lemma}\label{lem:classical}
Let $X$ be a derived stack with $\L_{X^\cl/X}\cong 0$. Then $i\colon X^\cl\to X$ is
an equivalence.
\end{lemma}

\begin{proof}
Being an equivalence is smooth local on $X$ and passing to the classical truncation commutes with smooth base
change, so we may assume $X=\Spec A$, and we must show $\pi_nA=0$ for all
$n\geqslant 1$. Suppose inductively that $\pi_kA=0$ for $1\leqslant k<n$. 
Then $\tau_{\geqslant 1}A=\tau_{\geqslant n}A$, so $\cofib(A\to\pi_0A)=\tau_{\geqslant n}A[1]$
is $(n+1)$-connective with $\pi_{n+1}=\pi_nA$. By Lemma~\ref{lem:conn} with
$m=n+1$,
\[
0=\pi_{n+1}\L_{\pi_0A/A}\cong \pi_0A\otimes_{\pi_0A}\pi_nA=\pi_nA .
\]
Hence $\pi_nA=0$ for all $n\geqslant 1$, i.e., the map $A\to\pi_0A$ is an equivalence.
\end{proof}

\begin{lemma}\label{lem:onedegree}
Let $M$ be a classical  stack. Let $C$ be a perfect complex on $M$ of
Tor-amplitude in $[a,b]$ whose cohomology vanishes in degrees $>a$. 
We have $C\cong\calE[-a]$ for a locally free sheaf $\calE$ of finite rank, and
$\chi(C)=(-1)^a\rk\calE$.
\end{lemma}

\begin{proof}
Since cohomology is concentrated in $[a,b]$ and vanishes
above $a$, the complex $C$ has a single cohomology sheaf $\calE=h^a(C)$.
It is coherent
because $C$ is perfect.
Thus $C\cong\calE[-a]$. For any discrete $N$ we have
\[
h^{a-j}\big(C\otimes^{\L}N\big)=\Tor_j(\calE,N),
\]
and the Tor-amplitude bound $\geqslant a$ forces $\Tor_j(\calE,N)=0$ for $j>0$. So
$\calE$ is flat and a flat coherent sheaf on a
scheme is locally free of finite rank. Finally
$\chi(\calE[-a])=(-1)^a\rk\calE$.
\end{proof}

\medskip\subsection{Proof of Theorem~\ref{thm:mainB}}

Let $d=\dim M=\vd X$.
Since $M$ is a smooth Artin stack of dimension $d$, the cotangent complex $\L_M$ is perfect of
Tor-amplitude in $[0,1]$ of rank 
$\chi(\L_M)=d$.
Consider the sequence of perfect complexes on $M$
\begin{equation}\label{eq:triangle}
E\xrightarrow{\ \phi\ }\L_M\longrightarrow \L_{M/X}.
\end{equation}
From \eqref{eq:triangle} we get a 
sequence 
$$\L_M\to \L_{M/X}\to E[1].$$ Here $\L_M$ has Tor-amplitude in $[0,1]$ and
$E[1]$ has Tor-amplitude in $[-2,0]$. The Tor-amplitude is preserved by extensions, so $\L_{M/X}$ is perfect of
Tor-amplitude in $[-2,1]$. By Lemma~\ref{lem:Lt0}, we have 
$$h^j(\L_{M/X})=0\quad\forall j>-2.$$
By Lemma~\ref{lem:onedegree} with $a=-2$, we have
$\L_{M/X}\cong\calE[2]$ with $\calE$ locally free, and the additivity of $\chi$ in
\eqref{eq:triangle} yields
\[
\rk\calE=\chi(\L_{M/X})=\chi(\L_M)-\chi(E)=d-\vd X=0 .
\]
A locally free sheaf of rank $0$ is zero, so $\L_{M/X}\cong0$ and
$\phi$ is an equivalence. This proves the part (a).

Since $\L_{M/X}\cong0$, Lemma~\ref{lem:classical} shows that
$i\colon M=X^\cl\to X$ is an equivalence. This proves the part (b).

The virtual structure sheaf of a quasi-smooth derived stack is
\[
[\calO_X^{\vir}]=\sum_{n\geqslant 0}(-1)^n[\pi_n\calO_X]\in K_0(X^\cl)
\]
By (b) we have $\pi_n\calO_X=0$ for
$n\geqslant 1$ and $\pi_0\calO_X=\calO_M$.
We deduce that $[\calO_X^{\vir}]=[\calO_M]$ in $K_0(M)$.
We now prove part (c) in Chow groups. First let $M$ be Deligne-Mumford. 
By \cite{BF},
the virtual class attached to the perfect obstruction theory
$\phi\colon E\to\L_M$ is
\[
[X]^{\vir}=0^{!}_{\mathfrak E}[\mathfrak C_M],\quad
\mathfrak E=h^1/h^0(E^\vee),
\]
where 
$$\mathfrak C_M\subset\mathfrak N_M=h^1/h^0(\T_M)$$ is the intrinsic normal cone
of $M$ and $\mathfrak C_M\hookrightarrow\mathfrak E$ is the closed embedding induced by
$\phi$. Since $M$ is smooth, we have
$$\mathfrak C_M=\mathfrak N_M=h^1/h^0(\T_M)\cong BT_M,$$
and by (a) the embedding $\mathfrak C_M\to\mathfrak E$ is the identity of $BT_M$.
Concretely, take the global resolution 
$$E^\vee\cong[E_0\to E_1]$$ with $E_0=T_M$ and $E_1=0$. 
Hence, the lifted cone $C_1\subset E_1=0$ vanishes, i.e., we have $C_1=M$, and
\[
[X]^{\vir}=0^{!}_{E_1}[C_1]=0^{!}_{M}[M]=[M].
\]
Equivalently, when $M$ is smooth and the obstruction sheaf
$\Ob=h^1(E^\vee)$ is locally free one has 
$$[X]^{\vir}=e(\Ob)\cap[M]$$
by \cite[\S5]{BF}. 
Here $E^\vee\cong\T_M$ by (a), so $\Ob=0$ and $e(\Ob)=1$.

For a general Artin stack, the virtual fundamental class of a quasi-smooth
derived stack is constructed in \cite{Khan} by deformation to the normal bundle
and Gysin pullback along the $1$-shifted cotangent bundle. By (b), the stack $X\cong M$ is
classical and smooth. For smooth $X$ this construction returns the
fundamental class $[X]$, see \cite{Khan}. Hence $[X]^{\vir}=[M]$.
\qed

\medskip

\begin{corollary}\label{cor:vir_generically_smooth}
Let $X$ be a quasi-smooth derived stack whose classical
truncation $M=X^\cl$ is quasicompact. Assume that
\begin{enumerate}[label=$\mathrm{(\alph*)}$,leftmargin=8mm,itemsep=1mm]
	\item $M$ is of pure dimension $\vd X$
	\item there is a dense open substack $M^\circ\subset M$ which is smooth.
\end{enumerate}
	Then $[X]^{\mathrm{vir}}=[M_{\mathrm{red}}]$ in $\H_{2\vd X}(M,\Q)$.
\end{corollary}

\begin{proof}
Set $d=\vd X$ and let $X^\circ\subset X$ be the open substack with
$(X^\cl)^\circ=M^\circ$. Then $X^\circ$ is quasi-smooth of virtual dimension $d$,
and by (i) and (ii) its truncation $M^\circ$ is smooth of dimension $d$.
By Theorem~\ref{thm:mainB} (c), $[X^\circ]^{\mathrm{vir}}=[M^\circ]$. Since virtual
fundamental classes are compatible with restriction to open substacks \cite{Khan}, we get
$$
[X]^{\mathrm{vir}}|_{M^\circ}=[M^\circ]=[M_{\mathrm{red}}]|_{M^\circ}.
$$
Let $Z=M\smallsetminus M^\circ$ with its reduced structure. By (i) and the density of $M^\circ$, we have $\dim Z\le d-1$, hence $\H_j(Z,\Q)=0$ for $j>2d-2$. The localization sequence
$$
\H_{2d}(Z,\Q)\to \H_{2d}(M,\Q)\to \H_{2d}(M^\circ,\Q)
\to \H_{2d-1}(Z,\Q)
$$
shows that restriction $\H_{2d}(M,\Q)\to \H_{2d}(M^\circ,\Q)$ is an isomorphism, and the claim follows.
\end{proof}

\medskip

\medskip\section{Tautological classes and the $\chi$-independence operator}\label{app:C}

We sketch a direct proof of \eqref{eq:tautperverspureHiggs2} based on \cite{MMSV} and \cite{DHKSV} . We begin with some recollection about the global BPS Lie algebra defined in \cite{DHKSV}, to which we refer for details. Consider the stack $\mathfrak{M}(r,d) \supseteq \mathfrak{M}^{\ss}(r,d)$ parametrizing \textit{all} Higgs sheaves on $C$ of rank $r$ and degree $d$. One may consider its Borel-Moore homology spaces $\H_i(\mathfrak{M}(r,d),\Q)$ and equip the direct sum 
$$\underline{\mathbf{A}}_C:=\bigoplus_{r,d} \H_{*+2(g-1)r^2}(\mathfrak{M}(r,d),\Q)$$ with the structure of an associative algebra structure. One also defines a compatible comultiplication 
$$\Delta: \underline{\mathbf{A}}_C \to \underline{\mathbf{A}}_C \hat{\otimes} \underline{\mathbf{A}}_C$$
and considers the Lie algebra of primitive elements\footnote{In \cite{DHKSV}, the global affinized BPS Lie algebra is denoted $\g_S$ and its semistable counterpart is denoted $\g_S^\ss$ where $S=T^*C$.}
$$\mathfrak{h}^C:=\bigoplus_{r,d} \mathfrak{h}^C(r,d), \qquad \mathfrak{h}^C(r,d)=\mathbf{A}_C(r,d)^{\text{prim}}.$$
The restriction map induces for every pair $(r,d)$ an isomorphism of vector spaces
$$\mathfrak{h}^C(r,d) \simeq \widehat{\g}^C(r,d) .$$
The $\chi$-independence operator is defined by
$$T_0([\text{pt}]):=\text{ad}\;D_{1,0}([\text{pt}]), \qquad D_{1,0}([\text{pt}])=[\text{pt}] \cap [\mathfrak{M}^\cl(0,1)] \in \mathfrak{h}_C(0,1).$$
It induces an isomorphism of vector spaces
$$T_0([\text{pt}]): \widehat{\g}^C(r,d) \simeq \mathfrak{h}^C(r,d) \stackrel{\sim}{\to} \mathfrak{h}^C(r,d+1) \simeq  \widehat{\g}^C(r,d+1) .$$
Furthermore, from the explicit description of the structure of the COHA of zero-dimensional sheaves (or a direct computation) it follows that for any $i \geqslant 1$
$$\left[T_0([\text{pt}]),\ch_i([\text{pt}])\right]=0$$
from which it also follows that $[T_0([\text{pt}]), \overline{c}_i([\text{pt}])]=0$ (note that this is false for general classes $\ch_i(\gamma)$).

\medskip

We now explain the proof of \eqref{eq:tautperverspureHiggs2}. Fix $(r,d)$ and first assume that $gcd(r,d)=1$. In that case we have $\widehat{\g}^C(r,d)=\H_*(\mathfrak{M}^\ss(r,d),\Q)$ and $\g^C(r,d)=\H_*(M^D(r,d),\Q)$ and \eqref{eq:tautperverspureHiggs2} follows from Poincar\'e duality together with the fact\footnote{In the notation of \cite{HMMS}, the passage from $c_k$ to $\overline{c}_k$ corresponds to the passage from $\psi_k$ to $Q_k$, see \textit{loc.cit.} Prop.~6.3.} proven in \cite{HMMS} that the perverse degree of $\overline{c}_k([\text{pt}])$ is equal to $k$. The case of a general pair $(r,d)$ is now a consequence of the facts that the $\chi$-independence operator commutes with $\overline{c}_k([\text{pt}])$ and that (for degree reasons) it maps $[\mathfrak{M}^{\ss,\rd}(r,d)]$ to a nonzero multiple of $[\mathfrak{M}^{\ss,\rd}(r,d+1)]$.

\end{document}